\documentclass[11pt]{amsart}

\usepackage{latexsym,amsmath,amsthm,amssymb,amscd,amsfonts,diagrams,graphics}

\DeclareFontFamily{U}{rsf}{}
\DeclareFontShape{U}{rsf}{m}{n}{
  <5> <6> rsfs5 <7> <8> <9> rsfs7 <10-> rsfs10}{}
\DeclareMathAlphabet{\mathscr}{U}{rsf}{m}{n}

\DeclareMathAlphabet{\mathgth}{U}{euf}{m}{n}

\DeclareFontFamily{U}{cyr}{}
\DeclareFontShape{U}{cyr}{m}{n}{
  <5> wncyr5 <6> wncyr6 <7> wncyr7 <8> wncyr8 <9> wncyr9 <10-> wncyr10}{}
\DeclareMathAlphabet{\mathcyr}{U}{cyr}{m}{n}

\input cyracc.def

\newcommand{\gCat}{\mathgth{Cat}}
\newcommand{\gGps}{\mathgth{Gps}}
\newcommand{\gVect}{\mathgth{Vect}}
\newcommand{\cA}{{\mathscr A}}

\newcommand{\sC}{{\mathcal C}}

\newcommand{\cC}{{\mathscr C}}
\newcommand{\OcC}{1\mbox{-}{\mathscr C}}
\newcommand{\cD}{{\mathscr D}}
\newcommand{\OcD}{1\mbox{-}{\mathscr D}}
\newcommand{\cE}{{\mathscr E}}
\newcommand{\cF}{{\mathscr F}}
\newcommand{\cG}{{\mathscr G}}
\newcommand{\cH}{{\mathscr H}}

\newcommand{\cO}{{\mathscr O}}

\newcommand{\cT}{{\mathscr T}}

\newcommand{\cX}{{\mathscr X}}
\newcommand{\cY}{{\mathscr Y}}

\newcommand{\FMYX}{\Phi_{Y\ra X}}

\newcommand{\FMXY}{\Phi_{X\ra Y}}

\newcommand{\FMYZ}{\Phi_{Y\ra Z}}
\newcommand{\FMXZ}{\Phi_{X\ra Z}}

\newcommand{\op}{{\mathrm{o}}}

\newcommand{\reg}{{\mathrm{reg}}}

\newcommand{\HKR}{{\mathrm{HKR}}}
\newcommand{\sHom}{\underline{\mathrm{Hom}}}

\newcommand{\D}{{\mathbf D}_{\mathrm{coh}}^b}
\newcommand{\PreD}{\tilde{\mathbf D}}

\newcommand{\chk}{{\scriptscriptstyle\vee}}
\newcommand{\R}{\mathbf{R}}
\newcommand{\Ld}{\mathbf{L}}
\newcommand{\lotimes}{\stackrel{\Ld}{\otimes}}

\newcommand{\pt}{{\mathrm{pt}}}

\newcommand{\bone}{{\mathbf 1}}
\DeclareMathOperator{\Tr}{Tr}
\DeclareMathOperator{\sTr}{Tr}
\DeclareMathOperator{\Nat}{Nat}
\DeclareMathOperator{\Spec}{Spec}
\DeclareMathOperator{\Rep}{Rep}

\DeclareMathOperator{\bExFun}{ExFun^+}
\DeclareMathOperator{\Prex}{Prex}

\newcommand{\Hom}{{\mathrm{Hom}}}
\newcommand{\THom}{\mbox{2-}{\mathrm{Hom}}}
\DeclareMathOperator{\End}{End}
\DeclareMathOperator{\Tor}{Tor}

\DeclareMathOperator{\Td}{td}
\DeclareMathOperator{\ch}{ch}

\DeclareMathOperator{\id}{id}

\DeclareMathOperator{\Id}{Id}

\DeclareMathOperator{\Ob}{Ob}

\DeclareMathOperator{\Ext}{Ext}

\newcommand{\adjoint}{\dashv}

\newcommand{\ra}{\rightarrow}
\newcommand{\lra}{\longrightarrow}

\newcommand{\scdot}{\,\cdot\,}

\newcommand{\C}{\mathbf{C}}

\newcommand{\Z}{\mathbf{Z}}

\newcommand{\gCoh}{\mathgth{Coh}}

\newcommand{\gQCoh}{\mathgth{Qcoh}}

\newcommand{\iso}{\cong}

\newcommand{\MP}[2]{\ensuremath{\langle #1,\, #2\rangle}}
\newarrow{Equal} =====

\theoremstyle{plain}
\newtheorem{theorem}{Theorem}[section]
\newtheorem{lemma}[theorem]{Lemma}
\newtheorem{corollary}[theorem]{Corollary}
\newtheorem{proposition}[theorem]{Proposition}

\theoremstyle{definition}
\newtheorem{definition}[theorem]{Definition}
\newtheorem{qdefinition}[theorem]{``Definition''}
\newtheorem{definition-theorem}[theorem]{Definition-Theorem}
\newtheorem{example}[theorem]{Example}

\theoremstyle{remark}
\newtheorem{remark}[theorem]{Remark}

\renewcommand{\phi}{\varphi}

\begin{document}

\author{Andrei C\u ald\u araru}

\title[The Mukai pairing, I]{The Mukai pairing, I: the Hochschild structure}

\date{}

\begin{abstract}
We study the Hochschild structure of a smooth space or orbifold,
emphasizing the importance of a pairing defined on Hochschild homology
which generalizes a similar pairing introduced by Mukai on the
cohomology of a K3 surface.  We discuss those properties of the
structure which can be derived without appealing to the
Hochschild-Kostant-Rosenberg isomorphism and Kontsevich formality,
namely:
\begin{itemize}
\item[--] functoriality of homology, commutation of push-forward with
the Chern character, and adjointness with respect to the generalized
pairing;
\item[--] formal Hirzebruch-Riemann-Roch and the Cardy condition from physics;
\item[--] invariance of the full Hochschild structure under
Fourier-Mukai transforms.
\end{itemize}
Connections with homotopy theory and TQFT's are discussed in an
appendix.  A separate paper~\cite{CalHH2} treats consequences of the
HKR isomorphism.  Applications of these results to the study of a
mirror symmetric analogue of Chen-Ruan's orbifold product will be
presented in a future paper.
\end{abstract}

\maketitle

\tableofcontents

\section{Introduction}
\label{sec:intro}

\subsection{}
The present work is the first in a series of three papers dedicated to
the study of the Hochschild structure of smooth spaces, laying out the
foundational material used in the other
two~\cite{CalHH2},~\cite{CalHH3}.  The Hochschild structure $(HH^*(X),
HH_*(X))$ is defined for a space $X$, and its fundamental properties
are studied.  The space $X$ can be an ordinary compact complex
manifold, or more generally a global quotient compact orbifold, a
proper Deligne-Mumford stack for which Serre duality holds, or a
compact ``twisted space'' in the sense of~\cite{Cal}.
\newpage

\subsection{}
\label{subsec:hochstr}
The {\em Hochschild structure} of $X$ consists of
\begin{itemize}
\item[--] a graded ring $HH^*(X)$, the Hochschild cohomology
  ring, whose $i$-th graded piece is defined as
\[ HH^i(X) = \Hom_{\D(X\times X)}(\cO_\Delta, \cO_\Delta[i]), \]
  where $\cO_\Delta=\Delta_*\cO_X$ is the structure sheaf of the
  diagonal in $X\times X$;
\item[--] a graded left $HH^*(X)$-module $HH_*(X)$, the Hochschild
  homology module, defined as
\[ HH_i(X) = \Hom_{\D(X\times X)}(\Delta_!\cO_X[i], \cO_\Delta), \]
  where $\Delta_!$ is the left adjoint of $\Delta^*$ defined by
  Grothendieck-Serre duality~(\ref{subsec:shriek});
\item[--] a non-degenerate graded pairing $\langle\scdot,\scdot\rangle$ on
  $HH_*(X)$, the {\em generalized Mukai pairing}.
\end{itemize}
The Hochschild cohomology ring has a rich and developed
theory~(\cite{GerSha},~\cite{Kon}).  The above definition of homology
is, to the author's knowledge, new (but see~\cite{WeiHH} for an
alternative equivalent definition, and~\cite{Mar} for a different
attempt).  The last important ingredient of the structure, the Mukai
pairing, has not been studied previously from the perspective of
Hochschild theory.

\subsection{}
In his groundbreaking work~\cite{MukK3} Mukai studied the relationship
between the derived category and the cohomology of K3 surfaces $X$ and
$Y$ by defining
\begin{itemize}
\item[--] a map $v:\D(X)\ra H^*(X, \C)$ given by
\[ v(\cF) = \ch(\cF).\Td(X)^{1/2}, \]
where $\Td(X)$ is the Todd genus of $X$ ($v(\cF)$ is called the {\em
Mukai vector} of $\cF$);
\item[--] an association $\Phi\mapsto\Phi_*$ which maps the integral
  transform 
\begin{align*}
 \Phi:\D(X)\ra \D(Y),\quad\quad & \Phi(\scdot) =
 \pi_{Y,*}(\pi_X^*(\scdot) \otimes \cE) 
\intertext{defined by an object $\cE$ in $\D(X\times Y)$ to the map on 
cohomology}
 \Phi_*:H^*(X, \C)\ra H^*(Y, \C), \quad\quad &\Phi_*(\scdot) =
 \pi_{Y,*}(\pi_X^*(\scdot) . v(\cE));
\end{align*}
\item[--] a pairing $\langle\scdot,\scdot\rangle$ on the cohomology
$H^*(X, \C)$, given by the formula
\[ \langle v, w\rangle = \int_X v_0. w_4 - v_2. w_2 + v_4.w_0, \]
where for a vector $v\in H^*(X,\C)$, $v_i$ is the component of $v$ in
$H^i(X, \C)$. 
\end{itemize}
It is worth emphasizing that the map $\Phi_*$ does not respect the usual
grading on the cohomology $H^*(X,\C)$.

\subsection{}
\label{subsec:Mukprop}
Mukai argued that the following properties are satisfied for K3
surfaces $X$ and $Y$:
\vspace*{2mm}
\begin{itemize}

\item[a.] {\bf Functoriality:} The association of maps on cohomology
  to integral transforms is functorial, in the sense that $\Id_{\D(X)}
  \mapsto \Id_{H^*(X, \C)}$, and $(\Phi\circ \Psi)_* = \Phi_*\circ \Psi_*$. 
\vspace{2mm}

\item[b.] {\bf Commutation with $v$:} The following diagram commutes
\[
\begin{diagram}[height=2em,width=2em,labelstyle=\scriptstyle]
\D(X) & \rTo^{\Phi} & \D(Y) \\
\dTo_v & & \dTo_v \\
H^*(X, \C) & \rTo^{\Phi_*} & H^*(Y, \C).
\end{diagram}
\]
\vspace*{2mm}

\item[c.] {\bf Adjointness:} 
If $\Phi:\D(Y)\ra \D(X)$ is left adjoint to $\Psi:\D(X)\ra \D(Y)$ then 
\[ \langle \Phi_* v, w \rangle_X = \langle v, \Psi_* w\rangle_Y \]
for $v\in H^*(Y, \C)$, $w\in H^*(X, \C)$.
\vspace*{2mm}

\item[d.] {\bf Hirzebruch-Riemann-Roch:} For $\cE, \cF\in \D(X)$ we
  have 
\[ \langle v(\cE), v(\cF) \rangle = \chi(\cE, \cF), \] 
where $\chi(\scdot,\scdot)$ is the Euler pairing on $K_0(X)$, 
\[ \chi(\cE, \cF) = \sum_i (-1)^i \dim \Ext^i_X(\cE, \cF). \]
\end{itemize}

It follows immediately from these properties that if $\Phi$ is an
equivalence of triangulated categories, then $\Phi_*$ is an isometry
between the corresponding cohomology groups, endowed with the Mukai
pairing. 

\subsection{}
This paper is devoted to generalizing Mukai's results to a wide class
of compact spaces, including in particular smooth compact complex
manifolds, twisted spaces in the sense of~\cite{Cal}, and certain
orbifolds or Deligne-Mumford stacks for which Serre duality holds.
The main point we want to emphasize is that the natural target for
defining Mukai's structure is not singular cohomology but rather
Hochschild homology.  Replacing singular cohomology by Hochschild
homology, we shall obtain all of Mukai's results for the wide class of
spaces above.

The first observation that hints to the fact that ordinary cohomology
is not the right target for the definition of the maps $\Phi_*$ is the
observation that in the case of a smooth compact complex manifold
these maps do not respect the usual grading on singular cohomology.
The correct grading that is preserved is the one given by the
verticals, and not the horizontals of the Hodge diamond of the space,
which is precisely the grading on Hochschild homology.

\subsection{}
\label{subsec:hkr}
To relate our approach to the original one of Mukai observe that the
Hochschild-Kostant-Rosenberg theorem asserts the existence of an
isomorphism~(\cite{CalHH2})
\[
\begin{diagram}[labelstyle=\scriptstyle]
I_\HKR:HH_i(X) \iso \bigoplus_{q-p=i} H^p(X, \Omega^q_X) &
 \rEqual^{\mathrm{\ def\ }} H\Omega_i(X) 
\end{diagram}
\]
between the $i$-th Hochschild homology of a smooth projective manifold
$X$ and the $n+i$-th column of the Hodge diamond of $X$.

It would seem natural to expect that, in the case of a K3 surface $X$,
the $I_\HKR$ isomorphism will match the abstract structures that we shall
define on $HH_*(X)$ with the original structure of Mukai. However, we
believe that a correction is needed for that: in~\cite{CalHH2} we
conjecture that we need to adjust the $I_\HKR$ isomorphism by
multiplying it by $\Td(X)^{1/2}$ before the abstract structure we define
will yield Mukai's original one.

\subsection{}
Let us now present our results.  After some generalities on integral
transforms and Serre duality in Section~\ref{sec:prel}, we discuss the
construction of left-right adjoint functors in
Section~\ref{sec:basiccon}.  This will be the basis for all our
results.  Then in Section~\ref{sec:defprops} we define Hochschild
homology and cohomology, as well as the generalized Mukai product.  In
Section~\ref{sec:functoriality} we consider an integral transform
$\Phi:\D(X) \ra \D(Y)$ between two spaces $X$ and $Y$ and define a
natural map of graded vector spaces $\Phi_*:HH_*(X) \ra HH_*(Y)$.
Using this construction we present in Section~\ref{sec:chernchar} a
definition of a Chern character map
\[ \ch: K_0(X) \ra HH_0(X) \]
which agrees, under $I_\HKR$, with the usual Chern character
map~(\cite{CalHH2})
\[ \ch: K_0(X) \ra \bigoplus_p H^p(X, \Omega^p_X). \]
(And which, under the corrected isomorphism $HH_0(X) \iso \bigoplus
H^{p,p}(X)$, yields the Mukai vector.  It is worth emphasizing that on
the level of Hochschild homology, no correction by the Todd genus is
needed; this correction appears in the usual statements because of the
``wrong'' choice of HKR isomorphism.)  We also discuss a definition
equivalent to ours given by Markarian~\cite{Mar}.

\subsection{}
The formal properties a --- d of~(\ref{subsec:Mukprop}) can now be
proven to hold in full generality, using $HH_*(X)$ instead of $H^*(X,
\C)$ and $\ch$ instead of $v$.  The corresponding results are
Theorems~\ref{thm:funct},~\ref{thm:commut},~\ref{thm:adjadj},
and~\ref{thm:hrr}.  A slightly more general version of the
Hirzebruch-Riemann-Roch theorem can be stated in this context
(Theorem~\ref{thm:cardy}).  Its origins can be traced to the Cardy
condition in physics.  It turns out in fact that properties a and c
are truly fundamental, while b and d are easy consequences of them.

\subsection{}
The final result of the paper is a proof, in
Section~\ref{sec:invderequiv}, of the fact that the full Hochschild
structure is invariant with respect to Fourier-Mukai transforms.  The
main result is:

\vspace{1mm}
\noindent
{\bf Theorem~\ref{thm:dcatinv}.}
\footnote{Parts of this theorem were also proven independently by
Orlov~\cite{Orlrus}.}  
{\em Let $X$ and $Y$ be spaces whose derived categories are equivalent
via a Fourier-Mukai transform (i.e., the equivalence is given by an
integral transform).  Then there exists a natural isomorphism of
Hochschild structures
\[ (HH^*(X), HH_*(X)) \iso (HH^*(Y), HH_*(Y)). \]
More precisely, there exists an isomorphism of graded rings
$HH^*(X)\iso HH^*(Y)$, an isomorphism $HH_*(X)\iso HH_*(Y)$ of graded
modules over the corresponding cohomology rings, and the latter
isomorphism is an isometry with respect to the Mukai pairings on
$HH_*(X)$ and on $HH_*(Y)$, respectively. }

\subsection{}
Throughout the paper there will be a certain tension between the
``$\Ext$'' interpretation of the Hochschild structure given
in~(\ref{subsec:hochstr}) and a parallel categorical interpretation.
The point is that there are alternative ways, outlined in
Appendices~\ref{app:categ1} and~\ref{app:categ2}, to regard elements
of $HH^*(X)$ and $HH_*(X)$ not as morphisms in $\D(X\times X)$, but
rather as natural transformations between certain functors $\D(X)\ra
\D(X)$.  Unfortunately, despite our best efforts, we have been unable
to make these ideas fully precise.  This appears to be primarily
caused by certain known technical problems with the definition of the
derived category~\cite{Gro}.  However the intuition behind the
categorical interpretation is most often the correct one, and as a
compromise we have decided to steer a middle course: we presented our
results in mathematically correct form in the ``$\Ext$''
interpretation, and we gave the intuitive ideas in the categorical
context.  We highly recommend the reader to read the appendices for
gaining intuition into the proofs.  

The current state of affairs is somewhat unsatisfactory, as several of
the proofs appear unnecessarily complicated.  We can only hope that
future developments of category theory will enable us to rewrite this
paper at a later date in the ``correct'' (categorical) language.

\vspace{1.5mm}
\noindent
\textbf{Acknowledgments.}  I have greatly benefited from conversations
with Eyal Markman, Tom Bridgeland, Mircea Musta\c t\u a, Andrew
Kresch, Tony Pantev, Jonathan Block, Justin Sawon and Sarah
Witherspoon.  Many of the ideas in this work were inspired by an
effort to decipher the little known but excellent work~\cite{Mar} of
Nikita Markarian.  Greg Moore suggested the connection between the
formal Riemann-Roch theorem and the Cardy condition in physics.  The
author's work has been supported by an NSF postdoctoral fellowship and
by travel grants and hospitality from the University of Pennsylvania,
the University of Salamanca, Spain, and the Newton Institute in
Cambridge, England.

\vspace{1.5mm}
\noindent
\textbf{Conventions.}  Throughout the paper a {\em space} will be a
compact complex manifold or proper algebraic variety over an
algebraically closed field of characteristic zero, possibly endowed
with an Azumaya algebra, or a smooth compact Deligne-Mumford stack
which satisfies Serre duality.  The derived category of a space will
refer to the bounded derived category of coherent sheaves on the
underlying space (which, in the case of the existence of an Azumaya
algebra $\cA$, shall mean coherent sheaves of modules over $\cA$).
Functors between derived categories shall always be implicitly
derived, but we shall keep clear the distinction between $\Hom$ and
$\R\Hom$.  Whenever we write $\cF\otimes \mu$ where $\cF$ is an object
and $\mu$ is a morphism, we mean $\Id_\cF\otimes \mu$.

\section{Preliminaries}
\label{sec:prel}

In this section we set up the basic context and notation.  We also
provide a brief introduction to Serre functors and Grothendieck-Serre
duality.  Our basic reference for these results is~\cite{BonKapSerre}.
We discuss a trace map that arises from the existence of Serre
functors and which is intimately related to one studied by
Illusie~\cite{Ill} and Artamkin~\cite{Art}.

\subsection{}
\label{subsec:inttransf}
Let $X$ and $Y$ be spaces, and let $\cE$ be an object in $\D(X\times
Y)$.  If $\pi_X$ and $\pi_Y$ are the projections from $X\times Y$ to
$X$ and $Y$, respectively, define the functor
\[ \FMXY^\cE:\D(X) \ra \D(Y) \quad \FMXY^\cE(\scdot) = \R\pi_{Y, *}
(\pi_X^*(\scdot) \lotimes \cE), \] 
which will be called the {\em integral functors} (on derived
categories) associated to $\cE$ (or with kernel $\cE$).

The association between objects of $\D(X\times Y)$ and integral
transforms is functorial: given a morphism $\mu:\cE\ra \cF$ between
objects of $\D(X\times Y)$, there is an obvious natural transformation
\[ \FMXY^\mu:\FMXY^\cE\Rightarrow \FMXY^\cF \]
given by
\[ \FMXY^\mu(\scdot) = \pi_{Y,*}(\pi_X^*(\scdot) \otimes \mu). \]

\subsection{}
\label{subsec:composite}
Given spaces $X, Y, Z$, and elements $\cE\in \D(X\times Y)$ and
$\cF\in \D(Y\times Z)$, define $\cF\circ\cE\in\D(X\times Z)$ by
\[ \cF\circ \cE =  \pi_{XZ, *} (\pi_{XY}^* \cE \otimes
\pi_{YZ}^* \cF), \] 
where $\pi_{XY}, \pi_{YZ}, \pi_{XZ}$ are the projections from $X\times
Y\times Z$ to the corresponding factors.  The reason behind the
notation is the fact that we have~(\cite[1.4]{BonOrl})
\[ \FMYZ^\cF \circ \FMXY^\cE = \FMXZ^{\cF\circ\cE}. \]

\subsection{}
\label{subsec:Serretrace}
Recall the definition of a (right) Serre functor on an additive
category $\sC$ with finite dimensional $\Hom$ spaces
from~\cite{ReiVan} (generalizing slightly the original definition
of~\cite{BonKapSerre}).  A Serre functor is a functor $S:\sC \ra \sC$
together with natural, bifunctorial isomorphisms
\[ \eta_{A,B}:\Hom_{\sC}(A, B) \stackrel{\sim}{\lra}
\Hom_{\sC}(B, SA)^\chk, \] 
for any $A,B$, where $\scdot^\chk$ denotes the dual vector space.  For
any $A$ in $\sC$, define
\[ \Tr:\Hom(A, SA) \ra \C,\quad\quad \Tr(f) = \eta_{A,A}(\id_A)(f). \]
\newpage

The following are easy consequences of the definition of a Serre
functor (see~\cite{ReiVan} for details):
\begin{lemma} 
For $f:A\ra B$ and $g:B\ra SA$, we have
\[ \eta_{A,B}(f)(g) = \Tr(g\circ f). \]
\end{lemma}

\begin{lemma}
\label{lem:shuffleok}
For $f:A\ra B$ and $g:B\ra SA$, we have
\[ \Tr(g\circ f) = \Tr(S f \circ g). \]
\end{lemma}

\subsection{}
To connect with more classical approaches to Serre duality, recall
that Illusie~(\cite{Ill}) and Artamkin~(\cite{Art}) construct a
trace map
\[ \Tr_\cE:\Hom_X(\cF\otimes\cE, \cG\otimes \cE) \ra \Hom_X(\cF, \cG)
\]
for objects $\cE, \cF, \cG$ in the derived category of a compact,
smooth space $X$, which generalizes the usual trace map on vector
spaces.  One way to write the definition of this trace map is that if
$\mu:\cF\otimes \cE \ra \cG\otimes \cE$ then $\Tr_\cE(\mu)$ is the
composition
\[
\begin{diagram}[height=2em,width=2em,labelstyle=\scriptstyle]
\cF & \rTo^{\id\otimes\eta} & \cF \otimes \cE \otimes \cE^\chk &
\rTo^{\mu\otimes \id} & \cG \otimes \cE \otimes \cE^\chk &
\rEqual^{\id \otimes \gamma} & \cG\otimes \cE^\chk\otimes \cE &
\rTo^{\id \otimes \epsilon} & \cG. 
\end{diagram}
\]
Here $\eta:\cO_X \ra \cE\otimes \cE^\chk \iso \R\sHom(\cE, \cE)$ is
the morphism which sends the section ``1'' of $\cO_X$ to the the
identity of $\Hom(\cE, \cE)$, $\gamma$ is the isomorphism that
interchanges the two factors, and $\epsilon$ is the original trace map
of Illusie and Artamkin.  This definition should be compared to the
generalized trace map of May~\cite{May}.

If we consider the functor
\[ S_X(-) = \omega_X[\dim X]\otimes\,-, \]
then in the standard form of Serre duality~(\cite{HarRD}) one constructs a
trace map
\[ \Tr_X:\Hom_X(\cO_X, S_X \cO_X) \ra \C \]
such that for objects $\cE$, $\cF$ of $\D(X)$ the pairing
\[ \langle \scdot, \scdot \rangle : \Hom_X(\cF, S_X \cE) 
\otimes \Hom_X(\cE, \cF) \ra \C \]
given by
\[ \langle f, g\rangle = \Tr_X(\Tr_\cE(f\circ g)) \]
is non-degenerate.  This yields isomorphisms
\[ \Hom_X(\cE, \cF) \iso \Hom_X(\cF, S_X\cE)^\chk \]
for every $\cE, \cF$ which are natural in both variables.
Applying~\cite[I.1.4]{ReiVan} it follows immediately that $S_X$ is a
Serre functor for $\D(X)$.  We shall often abuse notation and denote
$\Tr_X(\Tr_\cE(\scdot))$ by $\Tr_X(\scdot)$.

\subsection{}
The following generalization of a standard formula from linear algebra
is a rather involved result in category theory:
\begin{proposition}
\label{lem:tracetriangles}
Assume given a map of triangles
\[ 
\begin{diagram}[height=2em,width=2em,labelstyle=\scriptstyle]
\cE & \rTo & \cF & \rTo & \cG & \rTo & \cE[1] \\
\dTo^{e} & & \dTo^{f} & &
\dTo^{g} & & \dTo^{e[1]} \\
\cE\otimes\cH  & \rTo & \cF\otimes \cH & \rTo & \cG\otimes \cH & \rTo
& \cE\otimes \cH[1],
\end{diagram}
\]
where the bottom triangle is obtained by tensoring the top one with
$\cH$.  Assume furthermore that it is possible to find representatives
of all the objects in the diagram as complexes of sheaves, such that
the maps are maps of complexes and the squares commute on the nose
(not just up to homotopy), and that $g$ is the natural quotient map
obtained from $e$ and $f$.  Then we have
\[ \Tr_\cE(e) - \Tr_\cF(f) + \Tr_\cG(g) = 0 \]
as morphisms $\cO_X \ra \cH$.
\end{proposition}

\begin{proof}
This is~\cite[Theorem~1.9]{May}.  The condition about representing the
objects as complexes, etc., is precisely what the proof of [loc.\
cit.] uses.
\end{proof}

\noindent
The following is an easy exercise 
in linear algebra:
\begin{lemma}
\label{lem:tracetwotimes}
Let $\mu:\cE\otimes \cF\ra \cE\otimes \cG$ and $\nu:\cG\ra\cH$ be
morphisms in $\D(X)$.  Then
\[ \Tr_\cE((\id_\cE\otimes\, \nu) \circ \mu) = \nu \circ \Tr_\cE(\mu) \]
as morphisms $\cF\ra \cH$.
\end{lemma}

\section{The basic construction}
\label{sec:basiccon}

\subsection{}
The following rather innocuous remark about the construction of a
right adjoint functor from a left adjoint one is the basis for all the
results in this paper.  Consider a functor 
\[ \Phi:\D(X) \ra \D(Y) \]
that admits a left adjoint
\[\Phi^*:\D(Y)\ra \D(X). \]
Let $\cF$ and $\cG$ be objects in $\D(X)$.  The
fact that $\Phi$ is a functor implies that there is a natural map
\begin{align*}
\begin{diagram}[height=2em,width=3em,labelstyle=\scriptstyle]
\Hom_X(\cG, \cF) & \rTo^{\ \ \ \ \ \Phi\ \ \ \ \ \ \ \ } & \Hom_Y(\Phi\cG, \Phi\cF). 
\end{diagram}
\intertext{By Serre duality we can construct a left adjoint of this
  map (with
  respect to the Serre pairing)}
\begin{diagram}[height=2em,width=3em,labelstyle=\scriptstyle]
\Hom_X(\cF, S_X \cG) & \lTo^{\ \ \ \ \,\Phi^\dagger\ \ \ } & \Hom_Y(\Phi\cF, S_Y
\Phi\cG).
\end{diagram}
\end{align*}
The following proposition gives an explicit description of the map
$\Phi^\dagger$. 
\begin{proposition}
\label{prop:fundprop}
Let $\Phi^!:\D(Y) \ra \D(X)$ be given by 
\[ \Phi^! = S_X \circ \Phi^* \circ S_Y^{-1}. \]
Then $\Phi^!$ is a right adjoint to $\Phi$, and if
$\nu\in\Hom_Y(\Phi\cF, S_Y\Phi\cG)$ then $\Phi^\dagger\nu$ is the
composition
\[
\begin{diagram}[height=2em,width=2em,labelstyle=\scriptstyle]
\Phi^\dagger\nu:\cF & \rTo^{\bar{\eta}} & \Phi^!\Phi\cF & \rTo^{\Phi^!\nu}
& \Phi^! S_Y \Phi \cG & \rEqual & S_X \Phi^* \Phi\cG & \rTo^{S_X\epsilon}
& S_X \cG,
\end{diagram}
\]
where $\bar{\eta}$, $\epsilon$ are the unit and counit of the
adjunctions $\Phi\adjoint \Phi^!$, $\Phi^*\adjoint \Phi$,
respectively.

Explicitly, for $\mu\in\Hom_X(\cG, \cF)$ and $\nu\in\Hom_Y(\Phi\cF,
S_Y\Phi\cG)$ we have 
\[ \Tr_X(\Phi^\dagger \nu \circ \mu) = \Tr_Y(\nu \circ \Phi \mu). \]
\end{proposition}

\remark 
There is a striking similarity between the definition of
$\Phi^\dagger\nu$ and the definition in~\cite{May} of the generalized
trace maps.  It would be interesting to get a good explanation of this
similarity.

\begin{proof} 
Serre duality on $X$ and $Y$ gives the following diagram for 
$\cF\in \D(X)$, $\cH\in\D(Y)$
\[ 
\begin{diagram}[height=1.5em,width=3em,labelstyle=\scriptstyle]
\nu & \rMapsto & & & \Phi^! \nu \circ \bar{\eta}\\
\Hom_Y(\Phi\cF, \cH) & \rTo^{\sim} & \Hom_X(\cF, S_X\Phi^*S_Y^{-1}\cH)
& \rEqual & \Hom_X(\cF, \Phi^!\cH) \\
& & & & \\
\mbox{dual to} & & \mbox{dual to} & & \\ 
& & & & \\ 
\Hom_Y(S_Y^{-1}\cH, \Phi\cF) & \lTo^{\ \sim\ } &
\Hom_X(\Phi^*S_Y^{-1}\cH, \cF) \\ 
\Phi \bar{\mu} \circ \eta & \lMapsto & \bar{\mu}, & &
\end{diagram}
\]
where $\eta$ is the unit of $\Phi^*\adjoint\Phi$, and the top and bottom
rows are dual to each other and are given by the adjunctions
$\Phi\adjoint\Phi^!$, $\Phi^*\adjoint \Phi$.  It follows that
$S_X\Phi^*S_Y^{-1}$ is a right adjoint to $\Phi$ (see
also~\cite[Section 2]{Orl}).

Reading the duality between the top and bottom rows of the above diagram
we get
\[ \Tr_X(\Phi^!\nu\circ\bar{\eta}\circ \bar{\mu}) =
\Tr_Y(\nu\circ\Phi\bar{\mu}\circ \eta). \] 
If we take $\cH = S_Y\Phi\cG$, $\mu$ and $\nu$ as in the statement of
the proposition, and $\bar{\mu} = \mu\circ\epsilon$, then
$\Phi\bar{\mu}\circ\eta$ is nothing but $\Phi\mu$, and we conclude
that
\[ \Tr_X(\Phi^!\nu\circ\bar{\eta}\circ\mu\circ\epsilon) =
\Tr_Y(\nu\circ \Phi\mu). \]
By the commutativity property of the trace (Lemma~\ref{lem:shuffleok})
this can be rewritten as
\[ \Tr_Y(\nu\circ\Phi\mu) = \Tr_X(S_X\epsilon\circ
\Phi^!\nu\circ\bar{\eta}\circ\mu) = \Tr_X(\Phi^\dagger \nu \circ \mu). \]
\end{proof}

\subsection{}
If for some functor $\Psi:\D(X)\ra\D(X)$ we have a natural transformation 
\[ \nu:\Phi\Longrightarrow S_Y\Phi\Psi \]
then the above construction yields a new natural transformation
\[ 
\begin{diagram}[height=2em,width=2em,labelstyle=\scriptstyle]
\Phi^\dagger \nu: 1_X & \rImplies^{\bar{\eta}} & \Phi^!\Phi &
\rImplies^\nu & \Phi^!  S_Y\Phi\Psi & \rEqual & S_X \Phi^*\Phi\Psi &
\rImplies^\epsilon & S_X \Psi,
\end{diagram}
\]
and thus we can view $\Phi^\dagger$ as a map
\[ \Phi^\dagger:\Nat(\Phi, S_Y\Phi\Psi) \ra \Nat(1_X, S_X \Psi) \]
where $\Nat$ denotes the set of natural transformations between the
corresponding functors.  By Proposition~\ref{prop:fundprop} we have
for any $\mu:\Psi\cF\ra\cF$
\[ \Tr_X((\Phi^\dagger\nu)_\cF \circ \mu) = \Tr_Y(\nu_{\Phi\cF}
\circ \Phi\mu). \]

\subsection{}
\label{subsec:shriek}
The same kind of argument can be used to define a left adjoint to
$\Phi$ when a right adjoint is known.  For example,
\[ \Delta_! = S_{X\times X}^{-1} \Delta_* S_X \]
is a left adjoint to $\Delta^*$, where $\Delta:X\ra X\times X$ is the
diagonal embedding.

\subsection{}
\label{subsec:dualiso}
A similar kind of construction is the following: assume $\Phi$ and
$\Psi$ are functors from $\D(X)$ to $\D(Y)$, which admit right
adjoints $\Phi^!$, $\Psi^!$, respectively.  Then there exists a
natural isomorphism
\[ \tau:\Nat(\Phi, \Psi) \stackrel{\sim}{\lra} \Nat(\Psi^!, \Phi^!), \]
which maps $\mu: \Phi\Longrightarrow \Psi$ to the composite
\[
\begin{diagram}[height=1.5em,width=3em,labelstyle=\scriptstyle]
\Psi^! & \rImplies^{\bar{\eta}_\Phi} & \Phi^!\Phi\Psi^! & \rImplies^{\mu} &
\Phi^!\Psi\Psi^! & \rImplies^{\bar{\epsilon}_\Psi} & \Phi^!.
\end{diagram}
\]
Indeed, an inverse to $\tau$ is given by mapping $\nu:\Psi^!
\Longrightarrow \Phi^!$ to the natural transformation
\[
\begin{diagram}[height=1.5em,width=3em,labelstyle=\scriptstyle]
\Phi & \rImplies^{\bar{\eta}_\Psi} & \Phi\Psi^!\Psi & \rImplies^\nu &
\Phi\Phi^!\Psi & \rImplies^{\bar{\epsilon}_\Phi} & \Psi.
\end{diagram}
\]

\section{The Hochschild structure: definition and basic properties}
\label{sec:defprops}

\subsection{}
In this section we define Hochschild homology and cohomology for a
space.  The Mukai product is also introduced, together with its
categorical interpretation.  For simplicity of exposition we present
everything for a smooth compact scheme with no group action; the case
of an orbifold (or Deligne-Mumford stack) is obtained by thinking of
all the objects involved as equivariant.  For example the diagonal in
$X\times X$ will be viewed as an equivariant subvariety of $X\times
X$, all $\Ext$'s are computed in the category of equivariant sheaves,
etc.\ (see~\cite[Section 4]{BKR} for details).  We give some hints on
how to deal with general orbifolds in~(\ref{subsec:genorb}).
Similarly, the case of twisted spaces will be obtained by working in a
twisted derived category (with the observation that the diagonal can
also be viewed as an $(\alpha, \alpha)$-twisted sheaf, etc.), and
Serre functors make sense~\cite{Cal}.

\subsection{}
Let $X$ be a smooth, proper variety of dimension $n$ over $\C$.  The
following notations will be used throughout the paper:
\begin{itemize}
\item[--] $\Delta:X\ra X\times X$ is the diagonal embedding;
\item[--] $\omega_X$ is the canonical bundle of $X$;
\item[--] $S_X = \omega_X[n]$ as an object of $\D(X)$; often we shall
  also think of $S_X$ as the Serre functor $S_X \otimes\,-$\,;
\item[--] $S_{X\times X} = \omega_{X\times X}[2n]$ in $\D(X\times X)$; 
\item[--] $\cO_\Delta = \Delta_* \cO_X$, $S_\Delta = \Delta_* S_X$,
  $S_\Delta^{-1} = \Delta_* S_X^{-1}$;
\item[--] $\Delta_!:\D(X) \ra \D(X\times X)$ is the left adjoint of
  $\Delta^*$,
\[ \Delta_! = S_{X\times X}^{-1} \Delta_* S_X. \]
Note that $\Delta_!\cO_X \iso S_\Delta^{-1}$.
\end{itemize}

\begin{definition}
\label{def:hhext}
The Hochschild cohomology of $X$ is defined to be
\begin{align*}
HH^i(X) & = \Hom_{\D(X\times X)}(\cO_\Delta, \cO_\Delta[i]) =
\Ext^{i}_{X\times X} (\cO_\Delta, \cO_\Delta),
\intertext{and the Hochschild homology is defined as}
HH_i(X) & = \Hom_{\D(X)}(\Delta_! \cO_X[i], \cO_\Delta) =
\Ext^{-i}_{X\times X} (\Delta_!\cO_X, \cO_\Delta)
\\
& = \Ext^{-i}_{X\times X}(S_{X\times X}^{-1} \otimes S_\Delta,
\cO_\Delta) \\
& = \Ext^{-i}_{X\times X}(S_\Delta^{-1}, \cO_\Delta).
\end{align*}
\end{definition}

\subsection{}
This is a compact definition of the Hochschild groups, but for
completeness we include a discussion of the relationship of our
definition with the classical definitions of Weibel~\cite{WeiHH}.  For
further details on the cohomology side see~\cite{Swa}.

The bar resolution is defined to be the complex of quasi-coherent
sheaves of $\cO_{X\times X}$-modules
\[ \cdots \ra \cO_X^{\otimes n}\ra \cdots \ra \cO_X\otimes_\C \cO_X 
\otimes_\C \cO_X \ra \cO_X\otimes_\C\cO_X \ra 0, \]
with $\cO_{X\times X}$-module structure on $\cO_X^{\otimes n}$ given
by multiplication in the first and last factors, and with differential
\begin{eqnarray*} 
\lefteqn{d(a_0\otimes a_1\otimes\cdots\otimes a_n) = } \\
& & a_0a_1 \otimes a_2 \otimes\cdots \otimes a_n \,-\, 
a_0\otimes a_1a_2 \otimes \cdots \otimes a_n \,+\, \cdots \,+\, \\ 
& & (-1)^{n-1}a_0\otimes a_1\otimes\cdots \otimes
a_{n-1}a_n. 
\end{eqnarray*}
It is a resolution of $\cO_\Delta$ in $\gQCoh(X)$~\cite[1.1.12]{Lod}.

Hochschild cohomology is defined by Weibel by taking this resolution,
applying the functor $\sHom_{X\times X}(\scdot, \cO_\Delta)$, and then taking
hypercohomology of the resulting complex.  Since the bar resolution is
a resolution by free $\cO_X$-modules, applying $\sHom_{X\times
X}(\scdot, \cO_\Delta)$ and taking hypercohomology amounts to
computing the complex
\[ \R\Gamma(X\times X, \R\sHom_{X\times X}(\cO_\Delta, \cO_\Delta))
= \R\Hom_{X\times X}(\cO_\Delta, \cO_\Delta), \] 
whose $i$-th cohomology group is precisely $\Ext^{i}_{X\times
X}(\cO_\Delta, \cO_\Delta)$, which is our definition of $HH^i(X)$.  

Similarly, $HH_i(X)$ is usually defined by taking the bar resolution,
applying the functor $ - \otimes_{X\times X}\cO_\Delta$, and then
taking hypercohomology of the resulting complex, thought of as a
complex of $\cO_X$-modules by multiplication in the $\cO_\Delta$
factor.  (The complex obtained by tensoring the bar resolution with
$\cO_\Delta$ is usually referred to as the {\em bar complex}.)  In
derived category language this is equivalent to computing
\[ \R\Gamma(X,  \Delta^* \cO_\Delta) = \R\Hom_X(\cO_X, \Delta^*
\cO_\Delta) = \R\Hom_{X\times X}(\Delta_! \cO_X, \cO_\Delta). \]
Hence the $i$-th homology group of $\R\Gamma(X, \Delta^*
\cO_\Delta)$ (which is the classic definition of $HH_i(X)$) is
naturally isomorphic to the $i$-th homology (or $(-i)$-th cohomology)
group of
\[ \R\Hom_{X\times X}(\Delta_!\cO_X, \cO_\Delta). \]
which is our definition of $HH_i(X)$.  

\subsection{}
\label{subsec:genorb}
An alternative way of defining $HH_*(X)$ is to take the exact category
$\gCoh(X)$ of coherent sheaves on $X$ and to apply Keller's
construction~\cite{Kel}, which yields Hochschild homology.  This
provides an alternative easy way to define Hochschild homology for
arbitrary orbifolds: the usual notion of an orbibundle generalizes
immediately to that of a coherent orbisheaf, and these form an abelian
category.  Applying Keller's construction to this abelian category
yields a definition of Hochschild homology for an arbitrary orbifold.
A similar approach also works for the abelian category of twisted
sheaves.

\subsection{}
In the affine case, the idea of defining Hochschild homology as an
$\Ext$ group appears also in~\cite{Van} (where it is applied to the
study of Gorenstein rings, which are precisely the rings for which
Serre duality works as for smooth schemes).

\subsection{Degree Bounds}
The following result shows that homology and cohomology are non-zero
only in certain dimensions.
\begin{lemma}
If $\Delta$ is a locally complete intersection in $X\times X$ (in
particular, if $X$ is smooth), then
\[ \cH_i( \Delta^* \cO_\Delta) = 0 \]
for $i<0$ or $i>\dim X$, where $\cH_i(\Delta^*\cO_\Delta)$ denotes
the $i$-th homology sheaf of the complex $\Delta^*\cO_\Delta$.
\end{lemma}

\begin{proof}
The sheaf $\cH_i(\Delta^*\cO_\Delta)$ can be identified with
$\Tor_i^{X\times X}(\cO_\Delta, \cO_\Delta)$.  If $\cO_\Delta$ is a
locally complete intersection, this can be computed from the Koszul
resolution, which has length $\dim X$.  The result follows.
\end{proof}

From the lemma it follows immediately that cohomology lives in degrees
$0\leq i\leq 2n$ and homology can be non-zero only for $-n\leq i\leq
n$, where $n=\dim X$.  Indeed,
\[ \R\Hom_{X\times X}(\cO_\Delta, \cO_\Delta) \iso
\R\Hom_X(\Delta^* \cO_\Delta, \cO_X), \] 
and the Grothendieck spectral sequence computing the right-hand side
of the above equality will only have non-zero terms ${}^2E^{pq}$ in
the square $0\leq p,q\leq n$.  Similarly, the spectral sequence computing
\[ \R\Hom_X(\cO_X, \Delta^* \cO_\Delta) \]
(which yields Hochschild homology) will only have non-zero terms for
$0\leq p\leq n$, $-n\leq q\leq 0$.

\subsection{Ring-Module Structure}
Cohomology is naturally a graded ring, with product given by Yoneda
composition, and homology is a graded left $HH^*(X)$-module with the
same action.  The graded structure is given by the composition maps
\[ HH^i(X)\otimes HH_j(X) \ra HH_{j-i}(X). \]

\subsection{Mukai Product}
\label{subsec:mukprod}
Homology is equipped with a non-degenerate inner product (the
{\em Mukai product})
\[ \langle\scdot,\scdot\rangle:HH_*(X)\otimes HH_*(X) \ra \C, \]
which pairs $HH_i(X)$ with $HH_{-i}(X)$.  In order to define it, 
consider the contravariant functor
\[ {}^!:\D(X\times X) \ra \D(X\times X) \]
given by
\[ \cF \mapsto \cF^! = \rho(\cF^\chk \otimes \pi_1^* S_X), \]
where 
\[ \cF^\chk = \R\sHom_{X\times X}(\cF, \cO_{X\times X}). \]
and $\rho$ is the involution on $X\times X$ that interchanges the two
factors.  Since every object in $\D(X\times X)$ is quasi-isomorphic to
a finite complex of locally free sheaves of finite
rank,~\cite[II.5.16]{HarRD} shows that the functor ${}^!$ induces an
isomorphism
\[ \tau: \Hom_{X\times X}(\cF, \cG) \stackrel{\sim}{\lra} \Hom_{X\times X} 
(\cG^!, \cF^!), \]
for $\cF, \cG\in \D(X\times X)$.

\subsection{}
If we take $\cF = S_\Delta^{-1}[i]$ and $\cG = \cO_\Delta$, then
$\cG^! = \cO_\Delta$ and $\cF^! = S_\Delta[-i]$.  Indeed, we have 
\[ \cO_\Delta^\chk = \R\sHom_{X\times X}(\Delta_* \cO_X, \cO_{X\times
  X}) = \Delta_* \R\sHom_X(\cO_X, \Delta^! \cO_{X\times X}) = 
  \Delta_* S_X^{-1}, \]
where $\Delta^! = S_X \Delta^* S_{X\times X}^{-1}$ is the right adjoint of
$\Delta_*$~(\cite[III.11.1]{HarRD}), and thus
\[ \cO_\Delta^!  = \rho(\Delta_* S_X^{-1}\otimes \pi_1^* S_X) =
  \rho(\cO_\Delta) = \cO_\Delta, \] 
and similarly for $\cF^!$.  

Thus $\tau$ is an isomorphism between
\[ HH_i(X) = \Hom_{X\times X}(S_\Delta^{-1}[i], \cO_\Delta) \]
and $\Hom_{X\times X}(\cO_\Delta, S_\Delta[-i])$, which is the Serre
dual of $HH_{-i}(X)$,
\begin{align*}
\Hom_{X\times X}(\cO_\Delta, S_\Delta[-i]) & = \Hom_{X\times
  X}(\cO_\Delta, S_{X\times X}S_\Delta^{-1}[-i]) = 
\Hom_{X\times X}(S_\Delta^{-1}[-i], \cO_\Delta)^\chk \\
& = HH_{-i}(X)^\chk.
\end{align*}

\begin{definition}
\label{def:mukprod}
The non-degenerate pairing
\[ HH_i(X)\otimes HH_{-i}(X) \ra \C \]
given by 
\[ \langle v, w \rangle = \Tr_{X\times X}(\tau(v)\circ w). \]
is called the {\em generalized Mukai pairing}.  Note that it is not
symmetric in general.  
\end{definition}

\subsection{}
For a more intuitive (but not fully precise) introduction to the
Hochschild structure, the reader is suggested to consult
Appendices~\ref{app:categ1} and~\ref{app:categ2}.

\subsection{}
\label{subsec:tautensors}
In the particular case $\cF=S_\Delta^{-1}$, $\cG =
\cO_\Delta$,~(\ref{subsec:dualiso}) provides a better understanding of
the isomorphism
\[ \tau:\Hom_{X\times X}(S_\Delta^{-1}, \cO_\Delta) \stackrel{\sim}{\lra}
\Hom_{X\times X} (\cO_\Delta, S_\Delta). \] 
Indeed, $\tau$ will map $\mu:S_\Delta^{-1}[i]\ra \cO_\Delta$ to
\[
\begin{diagram}[height=1.5em,width=1.5em,labelstyle=\scriptstyle]
\tau(\mu):\cO_\Delta & \rTo^{\bar{\eta}} & S_\Delta[-i]
\circ S_\Delta^{-1}[i] \circ \cO_\Delta & \rTo^{\mu} & S_\Delta[-i] \circ
\cO_\Delta \circ \cO_\Delta & \rTo^{\bar{\epsilon}} & S_\Delta[-i], 
\end{diagram}
\]
where $\bar{\eta}$ and $\bar{\epsilon}$ are the ``unit'' and ``counit'' of the
``adjunctions'' $S_\Delta^{-1} \adjoint S_\Delta$, $\cO_\Delta\adjoint
\cO_\Delta$, respectively (see Proposition~\ref{prop:etaeps} for the
precise meaning of this unit and counit).  But we have
\[ S_\Delta[-i] \circ S_\Delta^{-1}[i] \circ \cO_\Delta \iso \cO_\Delta, \]
and thus $\bar{\eta}$ is a map $\cO_\Delta \ra \cO_\Delta$, which is
obviously the identity.  Similarly $\bar{\epsilon}$ is seen to be the
identity under the identification
\[ S_\Delta[-i]\circ \cO_\Delta\circ \cO_\Delta \iso S_\Delta[-i]. \]
We conclude that $\tau(\mu)$ is nothing else than $\mu\otimes
\pi_2^*S_X$, under the obvious identifications.  (We use $\pi_2$
because all the $S_\Delta$'s appear on the left.)

Observe that an identification similar to the one
in~(\ref{subsec:dualiso}) could be made using {\em left} adjoints
instead of right ones.  This gives another isomorphism
\[ \bar{\tau}:\Hom_{X\times X}(S_\Delta^{-1}[i], \cO_\Delta) \ra
\Hom_{X\times X}(\cO_\Delta, S_\Delta[-i]), \]
which is easily seen to be multiplication by $\pi_1^* S_X$.

\section{Functoriality of homology}
\label{sec:functoriality}

We present in this section the construction of a map of
graded vector spaces
\[ \Phi_*:HH_*(X) \ra HH_*(Y) \]
associated to an integral transform $\Phi:\D(X)\ra \D(Y)$. This
construction is natural in the sense that to the identity functor we
associate the identity map on homology, and $(\Phi\circ \Psi)_* =
\Phi_* \circ \Psi_*$ for composable integral transforms $\Phi$ and
$\Psi$.  It is worth pointing out that, despite its name, Hochschild
cohomology is not functorial in any reasonable sense.

\subsection{}
\label{subsec:pullback}
Let $\Phi:\D(X)\ra \D(Y)$ be an exact functor which admits a left
adjoint (for example, any integral transform).  Given an element
$\mu\in HH_*(X)$ we want to define $\Phi_*\mu$ in a way that would be
natural with respect to $\Phi$.

Let us begin with the categorical interpretation, where things are
easier.  Recall that in Section~\ref{sec:basiccon} we constructed a
natural map
\[ \Phi^\dagger:\Nat(\Phi, S_Y\Phi\Psi) \ra \Nat(1_X, S_X \Psi). \]
If we take $\Psi$ to be the shift functor $[i]$, there is a natural
restriction map
\[ \Nat(1_Y, S_Y[i]) \ra \Nat(\Phi, S_Y\Phi[i]), \]
and composing we get a map
\[ \Phi^\dagger:\Nat(1_Y, S_Y[i]) \ra \Nat(1_X, S_X[i]). \]
The defining property of $\Phi^\dagger$ is the equality
\[ \Tr_X((\Phi^\dagger \nu)_\cF \circ \mu) = \Tr_Y(\nu_{\Phi\cF} \circ
\Phi\mu) \] 
for any $\nu\in\Nat(1_Y, S_Y[i])$ and $\mu:\cF[i] \ra \cF$.  In
particular, if $\nu\in\Nat(1_Y, S_Y)$ and we take $\mu=\id_\cF$, we
have
\[ \Tr_X((\Phi^\dagger\nu)_\cF) = \Tr_Y(\nu_{\Phi\cF}). \]

To construct the map $\Phi_*:HH_i(X) \ra HH_i(Y)$ we would want to
take the adjoint of $\Phi^\dagger$ (recall that homology is thought of
as the dual of $\Nat(1_Y, S_Y)$ with respect to Serre duality of
natural transformations).  Unfortunately we do not know how to make
this precise, hence we need to switch to the ``$\Ext$''
interpretation.

\subsection{}
We want to use Proposition~\ref{prop:fundprop} to rewrite the above
definition in a way that generalizes to the ``$\Ext$'' interpretation.
Indeed, we want to find a map 
\[ \Phi_* : \Hom_{X\times X}(S_{\Delta_X}^{-1}[i], \cO_{\Delta_X}) \ra 
\Hom_{Y \times Y}(S_{\Delta_Y}^{-1}[i], \cO_{\Delta_Y}), \]
and not just a map on natural transformations.

By Proposition~\ref{prop:fundprop} we see that for $\nu\in\Nat(1_Y,
S_Y[i])$, $\Phi^\dagger\nu$ can be written as the composite
\[ 
\begin{diagram}[height=2em,width=1.5em,labelstyle=\scriptstyle]
\Phi^\dagger \nu: 1_X & \rImplies^{\bar{\eta}} & \Phi^!\Phi &
\rImplies^\nu & \Phi^!  S_Y\Phi[i] & \rEqual & S_X \Phi^*\Phi[i] &
\rImplies^\epsilon & S_X[i],
\end{diagram}
\] 
where $\bar{\eta}$ and $\epsilon$ are the unit and counit of the
respective adjunctions.  

Assume that $\Phi$ is an integral transform, given by an object
$\cF\in\D(X\times Y)$, and define
\begin{align*}
\cG = \cF^\chk\otimes \pi_Y^* S_Y, & \quad \cH = \cF^\chk\otimes
  \pi_X^* S_X.
\intertext{Then by~\cite[Lemma 4.5]{BriFM},}
\Phi^* = \FMYX^{\cG},&\quad \Phi^!=\FMYX^{\cH}
\end{align*}
are left and right adjoints of $\Phi$, respectively.  

\begin{proposition}
\label{prop:etaeps}
There exist natural morphisms $\bar{\eta}:\cO_{\Delta_X} \ra \cH\circ \cF$ and
$\epsilon:\cF\circ \cG\ra\cO_{\Delta_Y}$ that correspond to 
\[ \bar{\eta}:1_X\Rightarrow \Phi^!\circ\Phi\quad\mbox{ and }\quad\epsilon:
\Phi^*\circ\Phi\Rightarrow 1_Y \] 
under the correspondence between morphisms between objects on a
product and natural transformations of the underlying functors.
\end{proposition}

\begin{proof}
Let $\pi_{ij}$ be the projection from $X\times Y\times X$ onto the
$i$-th and $j$-th factors, so that
\[ \cH\circ \cF = \pi_{13,*}(\pi_{12}^*\cF \otimes \pi_{23}^* \cH). \]
Also, let 
\[ \Delta:X\times Y\ra X\times Y\times X \]
be the map that on points is given by
\[ (x,y) \mapsto (x,y,x). \]
Then we have
\begin{align*}
\Hom_{X\times X}(\cO_{\Delta_X}, \cH\circ \cF) & = \Hom_{X\times
  X}(\cO_{\Delta_X}, \pi_{13,*}(\pi_{12}^*\cF \otimes \pi_{23}^* \cH))
  \\
& = \Hom_{X\times Y\times X}(\pi_{13}^* \cO_{\Delta_X}, \pi_{12}^*\cF
  \otimes \pi_{23}^* \cH) \\
& = \Hom_{X\times Y\times X}(\Delta_*\cO_{X\times Y}, \pi_{12}^*\cF
  \otimes \pi_{23}^* \cH) \\
& = \Hom_{X\times Y}(\cO_{X\times Y}, \Delta^!(\pi_{12}^* \cF \otimes
  \pi_{23}^*(\cF^\chk \otimes \pi_X^*S_X))) \\
& = \Hom_{X\times Y}(\cO_{X\times Y}, \cF\otimes \cF^\chk) \\
& = \Hom_{X\times Y}(\cF, \cF),
\end{align*}
and we take $\bar{\eta}$ to be the image of the identity morphism of
$\cF$ under the above isomorphism.

The construction of $\epsilon$ is entirely similar and will be left to
the reader.
\end{proof}

\begin{definition}
Given a morphism $\nu:\cO_{\Delta_Y}\ra S_{\Delta_Y}[i]$ define
$\Phi^\dagger\nu:\cO_{\Delta_X}\ra S_{\Delta_X}[i]$ to be the
composite morphism
\[
\begin{diagram}[height=2em,width=2em,labelstyle=\scriptstyle]
\cO_{\Delta_X} & \rTo^{\bar{\eta}} & \cH\circ \cF=\cH\circ
\cO_{\Delta_Y}\circ \cF & \rTo^\nu &
\cH\circ S_{\Delta_Y}[i] \circ \cF = S_{\Delta_X}[i]\circ \cG\circ \cF & 
\rTo^\epsilon & S_{\Delta_X}[i],
\end{diagram}
\]
where $\bar{\eta}$ and $\epsilon$ are the maps defined in
Proposition~\ref{prop:etaeps}.  Define
\[ \Phi_*:HH_*(X)\ra HH_*(Y) \] 
as the right adjoint to the map $\Phi^\dagger$ with respect to Serre duality
on $X\times X$ and on $Y\times Y$, i.e., for $\mu\in HH_*(X)$,
$\Phi_*\mu$ is the unique element in $HH_*(Y)$ such that
\[ \Tr_{X\times X}(\Phi^\dagger\nu \circ \mu) = \Tr_{Y\times Y}(\nu\circ
\Phi_* \mu) \]
for every $\nu\in\Hom^*_{Y\times Y}(\cO_{\Delta_Y}, S_{\Delta_Y})$. 
\end{definition}

\subsection{}
The following theorem summarizes the functoriality properties of this
construction:
\begin{theorem}
\label{thm:funct}
The map on homology associated to the identity functor is the
identity, and if $\Psi:\D(X) \ra \D(Y)$ and $\Phi:\D(Y)\ra \D(Z)$ are
integral transforms then we have
\[ (\Phi\circ\Psi)_*=\Phi_*\circ\Psi_*. \]
\end{theorem}

\begin{proof}
Follows easily from the observation that if $\Phi^*$, $\Psi^*$ are
left adjoints to $\Phi$ and $\Psi$, then $\Psi^*\circ \Phi^*$ is a
left adjoint to $\Phi\circ \Psi$, and similarly for right adjoint.
Also, the obvious relations between units and counits hold.  This
proves the result at a categorical level, and we leave to the patient
reader the task of checking that the corresponding compatibilities
hold in the Ext interpretation.
\end{proof}

\section{The Chern character}
\label{sec:chernchar}

In this section we define the Chern character map $\ch:K_0(X) \ra
HH_0(X)$.  We also discuss an equivalent construction of
Markarian~\cite[Definition 2]{Mar}.

\subsection{}
Let $p:X\ra \pt$ be the structure map of $X$, and let $\cF$ be any
object in $\D(X)$.  Consider the functor $\Phi=\Phi^\cF_{\pt\ra X}$
defined by $\cF$, and observe that we have $\Phi(\cO_\pt) = \cF$.  By
the results in Section~\ref{sec:functoriality}, the integral transform
$\Phi$ induces a map on homology
\[ \Phi_*:HH_*(\pt) \ra HH_*(X). \]
Observe that $HH_0(\pt) = \Hom_{\pt\times \pt}(\cO_\pt, \cO_\pt)$ has
a distinguished element $\bone$ given by the identity (we use the fact
that $S_\pt = \cO_\pt$).

\begin{definition}
Define the Chern character of $\cF$, $\ch(\cF)$, by 
\[ \ch(\cF) = (\Phi^\cF_{\pt \ra X})_*(\bone) \in HH_0(X). \]
\end{definition}

\subsection{}
Since this definition is slightly hard to work with, we unravel it to
a more usable version.  Recall that in~(\ref{subsec:pullback}) we
defined a map 
\[ \Phi^\dagger:\Hom_{X\times X}(\cO_{\Delta_X},
S_{\Delta_X}) \ra \Hom_{\pt\times \pt}(\cO_{\Delta_\pt},
S_{\Delta_\pt}). \] 
Applying its defining property with $\mu = \id_{\cO_\pt}$, we get
\[ \Tr_{\pt\times\pt}(\Phi^\dagger \nu) = \Tr_X(\nu_{\Phi \cO_\pt}) =
\Tr_X(\nu_\cF) \] 
for any $\nu\in\Hom_{X\times X}(\cO_\Delta, S_\Delta)$. (Here
$\nu_\cF:\cF \ra S_X \cF$ is the value at $\cF$ of the natural
transformation induced by $\nu$.)  The map $\Phi_*$ that we are
interested in is the adjoint of $\Phi^\dagger$ with respect to
Serre duality on $X\times X$ and $\pt\times \pt$, respectively.
Explicitly, we must have the equality
\[ \Tr_{X\times X}(\nu\circ \Phi_*\bone) = \Tr_{\pt\times
  \pt}(\Phi^\dagger \nu\circ \bone) \]
and thus since $\bone$ is nothing but the identity, we conclude that
we must have
\[ \Tr_{X\times X}(\ch(\cF) \circ \nu) = \Tr_{\pt\times
  \pt}(\Phi^\dagger \nu) = \Tr_X(\nu_\cF). \]

\subsection{}
\label{subsec:iotae}
We rewrite the above definition as follows: a homomorphism
$\nu:\cO_\Delta \ra S_\Delta$ in $\D(X\times X)$ induces a natural
transformation
\[ \iota(\nu):1_X \Longrightarrow S_X \]
between the identity functor and the Serre functor on $\D(X)$.  Thus
for every $\cF\in\D(X)$ we get a map
\begin{align*}
\begin{diagram}[height=2em,width=2em]
\Hom_{X\times X}(\cO_\Delta, S_\Delta) & \rTo^{\ \ \iota_\cF\ \ } & \Hom_X(\cF,
S_X \cF),
\end{diagram}
\intertext{whose left adjoint with respect to the Serre duality pairing we 
denote by $\iota^\cF$:} 
\begin{diagram}[height=2em,width=2em]
HH_0(X) = \Hom_{X\times X}(S_\Delta\otimes S_{X\times X}^{-1},
\cO_\Delta) & \lTo^{\ \ \ \ \iota^\cF\ \ \ \ \ } & \Hom_X(\cF, \cF).
\end{diagram}
\end{align*}

With these notations, the above calculations reduce to the following
equivalent definition of $\ch(\cF)$, similar to one given by
Markarian~\cite{Mar}:
\begin{definition}
\label{def:chernchar}
The Chern character of $\cF$ is defined as the image
\[ \ch(\cF) = \iota^\cF(\id_\cF) \in HH_0(X) \]
of the identity morphism of $\cF$ in $HH_0(X)$ under $\iota^\cF$.  
Explicitly, $\ch(\cF)$ is the unique element of $HH_0(X)$ such that
\[ \Tr_{X\times X} (\nu \circ \ch(\cF)) = \Tr_X(\iota_\cF(\nu)) =
\Tr_X(\pi_{2,*}(\pi_1^*\cF \otimes \nu)) \] 
for all $\nu\in\Hom_{\D(X\times X)}(\cO_\Delta, S_\Delta)$.
\end{definition}

Under the Hochschild-Kostant-Rosenberg isomorphism this definition of
$\ch(\cF)$ agrees with the usual one~\cite[Theorem~4.5]{CalHH2}.

\subsection{}
The following proposition shows that the map $\ch:\D(X) \ra HH_0(X)$
factors through $\D(X)\ra K_0(X)$ to yield the desired Chern character
map
\[ \ch:K_0(X) \ra HH_0(X). \]

\begin{proposition}
If $\cF\ra \cG\ra \cH\ra \cF[1]$ is an exact triangle in $\D(X)$, then 
\[ \ch(\cF) - \ch(\cG) + \ch(\cH) = 0 \]
in $HH_0(X)$.
\end{proposition}

\begin{proof}
For any $\nu\in \Hom_{X\times X}(\cO_\Delta, S_\Delta)$, $\iota(\nu)$ is a
natural transformation, and as such it gives a map of triangles
\[ 
\begin{diagram}[height=2em,width=2em]
\cF & \rTo & \cG & \rTo & \cH & \rTo & \cF[1] \\
\dTo_{\iota_\cF(\nu)} & & \dTo_{\iota_\cG(\nu)} & &
\dTo_{\iota_\cH(\nu)} & & \dTo_{\iota_\cF(\nu)[1]} \\
S_X \cF & \rTo & S_X \cG & \rTo & S_X \cH & \rTo & S_X \cF[1].
\end{diagram}
\]
Observe that if we represent $\nu$ by an actual map of complexes of
injectives, and $\cF$, $\cG$, $\cH$ by complexes of locally free
sheaves, then the resulting maps in the above diagram commute on the
nose (no further injective or locally free resolutions are needed), so
we can apply Proposition~\ref{lem:tracetriangles} to get
\[ \Tr_X(\iota_\cF(\nu)) - \Tr_X(\iota_\cG(\nu)) + \Tr_X(\iota_\cH(\nu))
= 0. \]
Therefore
\[ \Tr_{X\times X}(\nu\circ(\ch(\cF)-\ch(\cG)+\ch(\cH))) = 0 \]
for any $\nu\in \Hom_{X\times X}(\cO_\Delta, S_\Delta)$; since the
Serre duality pairing between $\Hom_{X\times X}(\cO_\Delta, S_\Delta)$
and $HH_0(X)$ is non-degenerate, we conclude that
\[ \ch(\cF) - \ch(\cG) + \ch(\cH) = 0. \]
\end{proof}

\begin{example} 
To have a non-commutative example at hand, consider the case when $G$
is a finite group, acting trivially on a point.  The resulting
orbifold $BG = [\scdot/G]$ can be thought of as $\Spec R$ where
$R=\C[G]$, the group ring of $G$.  Indeed, a coherent sheaf on $BG$,
which by definition is a finite-dimensional representation of $G$, is
precisely the same thing as a module over $\C[G]$.  The Serre functor
on $BG$ is trivial, and $BG\times BG$ should be thought of as $\Spec
(R\otimes R^\op)$, with $\cO_\Delta$ represented by $R$ as a module
over $R\otimes R^\op$ by left and right multiplication ($R^\op$
denotes the opposite ring of $R$).  Thus
\[ HH_0(BG) = \Hom_{R\otimes R^\op}(R, R) = Z(R), \]
where $Z(R)$ represents the center of $R$ (composition of morphisms in
$\Hom(R, R)$ is the same as multiplication in $Z(R)$ under the
identification).  Thus the Chern character is a map from
$K_0(\Rep(G))$ to the center $Z(R)$ of the group ring.

To understand this map, let $f\in Z(R)$ and $V$ be a representation of
$G$ (i.e., a right $R$-module).  The map $\iota^V\!(f):V\ra V$ is
multiplication by $f$ on the right.  The Chern character of
$V$, $\ch(V)$, is by definition the unique element $e_V\in Z(R)$ such
that 
\[ \Tr_{BG\times BG}(e_V\cdot f) = \Tr_{BG}(\iota^V\!(f)) \]
for all $f\in Z(R)$.  The left hand side is $\chi_{V_\reg}(e_V \cdot
f)$, the value at $e_V\cdot f$ of the character $\chi_{V_\reg}$ of the
regular representation $V_\reg=R$ of $G$, and the right hand side is
the value of $\chi_V$ at $f$.  Recall that $R$, being semisimple, is
isomorphic to the direct sum of the endomorphism algebras $\End(V_i)$
over a set of representatives $\{V_i\}$ of isomorphism classes of
irreducible representations of $G$ (Wedderburn's theorem).  Let
$\{e_i\}$ be the orthogonal set of idempotents corresponding to this
decomposition.  Then it is obvious from the fact that multiplication
by $e_i$ is the projection on the $\End(V_i)$ component that we have
\[ \chi_{V_\reg}(e_i\cdot f) = \chi_{V_i}(f), \]
for any $f\in Z(R)$, and thus it follows that 
\[ \ch(V_i) = e_i. \]
By semisimplicity this computes the value of the Chern character of
any representation.

The explicit value of $\ch(V_i)$ can be found in~\cite[2.12]{Isa}:
\[ \ch(V_i) = \frac{1}{|G|} \sum_{g\in G} \chi_{V_i}(1)\chi_{V_i}(g^{-1})
g. \]
\end{example}

\section{Properties of the structure}

In this section we argue that properties b, c and d of the original Mukai
construction hold if we replace $H^*(X,\C)$ with Hochschild homology
and  $v$ by $\ch$. 

\subsection{}
The commutativity of $\Phi_*$ and $\ch$ is the content of the
following theorem:
\begin{theorem}
\label{thm:commut}
The following diagram commutes for any integral transform $\Phi:\D(X)\ra\D(Y)$:
\[
\begin{diagram}[height=2em,width=2em,labelstyle=\scriptstyle]
\D(X) & \rTo^{\Phi} & \D(Y) \\
\dTo_{\ch} & & \dTo_{\ch} \\
HH_0(X) & \rTo^{\Phi_*} & HH_0(Y).
\end{diagram}
\]
\end{theorem}

\begin{proof}
We use the first of the two equivalent definitions of $\ch$ given in
Section~\ref{sec:chernchar}.  Let $\cE$ be an object of $\D(X)$, and
let $\cF = \Phi\cE$.  Observe that we have 
\[ \Phi\circ \Phi_{\pt\ra X}^\cE = \Phi_{\pt \ra Y}^\cF, \]
since any functor from $\D(\pt)$ is determined by its value at
$\cO_\pt$.  By Theorem~\ref{thm:funct} we have
\begin{align*}
\Phi_* \ch(\cE) & = \Phi_* [(\Phi_{\pt \ra X}^\cE)_* \bone] =
(\Phi\circ \Phi_{\pt\ra X}^\cE)_* \bone \\ 
& = (\Phi_{\pt \ra Y}^\cF)_* \bone = \ch(\cF) = \ch(\Phi\cE).
\end{align*}
\end{proof}
%% Let $\cF$ be an object in $\D(X)$.  Then by the definition of the
%% Chern character~(\ref{def:chernchar}), $\ch(\Phi\cF)$ is the unique
%% element in $HH_0(Y)$ such that
%% \[ \Tr_{Y\times Y}(\mu\circ \ch(\Phi\cF)) =
%% \Tr_Y(\iota_{\Phi\cF}(\mu)) \] 
%% for any $\mu:\cO_{\Delta_Y} \ra S_{\Delta_Y}$.  If we let $\nu =
%% \iota(\mu)$, the main property of $\Phi^\dagger$
%% in~(\ref{subsec:pullback}) can be read as
%% \[ \Tr_Y(\iota_{\Phi\cF}(\mu)) = \Tr_Y(\nu_{\Phi\cF}) =
%% \Tr_X((\Phi^\dagger\nu)_\cF) = \Tr_X(\iota_\cF(\Phi^\dagger\mu)).\]
%% (We make use here of the compatibility between the ``Ext'' and the
%% natural transformations interpretations.)  We conclude that for any
%% $\mu:\cO_{\Delta_Y} \ra S_{\Delta_Y}$ we have
%% \begin{align*}
%% \Tr_{Y\times Y}(\mu\circ \ch(\Phi\cF)) & = \Tr_Y(\iota_{\Phi\cF}(\mu))
%% \\
%% & = \Tr_X(\iota_\cF(\Phi^\dagger\mu)) \\
%% & = \Tr_{X\times X}(\Phi^\dagger\mu \circ \ch(\cF)) \\
%% & = \Tr_{X\times X}(\mu \circ \Phi_*\ch(\cF)). 
%% \end{align*}
%% Since $\ch(\Phi\cF)$ is uniquely determined by its dual functional, it
%% follows that 
%% \[ \ch(\Phi\cF) = \Phi_*\ch(\cF). \]

\subsection{}
We now move on to adjoint properties of maps on homology induced by
adjoint functors. 
\begin{proposition}
\label{prop:xytoxx}
Let $\mu:\cF\ra S_{X\times Y}\cF = S_{\Delta_Y}\circ \cF \circ
S_{\Delta_X}$, and let $\mu'$ be the composite morphism
\begin{align*}
&\begin{diagram}[height=2em,width=2em,labelstyle=\scriptstyle]
\mu':\cO_{\Delta_X} & \rTo^{\bar{\eta}} & \cH\circ \cF &
\rTo^{\cH\circ\mu} & \cH\circ S_{\Delta_Y} \circ \cF \circ
S_{\Delta_X} \iso S_{\Delta_X}\circ \cG\circ \cF\circ S_{\Delta_X} &
\rTo^{S_{\Delta_X}\circ\epsilon\circ S_{\Delta_X}} 
\end{diagram} \\
& \quad \rTo  S_{\Delta_X}\circ
S_{\Delta_X} = S_{X\times X} \cO_{\Delta_X}.
\end{align*}
Then 
\[ \Tr_{X\times Y}(\mu) = \Tr_{X\times X}(\mu'). \]
\end{proposition}

\begin{proof}
Follows from a calculation entirely similar to that of
Proposition~\ref{prop:fundprop} which is left to the
reader.
\end{proof}

\begin{theorem}
\label{thm:adjadj}
Let $\Psi:\D(Y)\ra \D(X)$ be a left adjoint to $\Phi$.  Then we have
\[ \langle v, \Phi_* w \rangle = \langle \Psi_* v, w\rangle \]
for $v\in HH_*(Y)$, $w\in HH_*(X)$.
\end{theorem}

\begin{proof}
We begin with the observation that it is enough to show that 
\[ \tau \Psi_* = \Phi^\dagger \tau:HH_*(Y) \ra \Hom_{X\times
  X}^*(\cO_{\Delta_X}, S_{\Delta_X}), \]
where $\tau$ is the map defined in~(\ref{subsec:mukprod}).  Indeed, if
this equality holds, we have 
\begin{align*}
\langle v, \Phi_* w \rangle & = \Tr_{Y\times Y}(\tau(v)\circ
\Phi_* w) = \Tr_{X\times X}(\Phi^\dagger(\tau(v))\circ w) \\
& = \Tr_{X\times X}(\tau(\Psi_*(v))\circ w) = \langle \Psi_*v,
w\rangle. 
\end{align*}

We have observed in~(\ref{subsec:tautensors}) that
\[ \tau:\Hom_{X\times X}(S_\Delta^{-1}, \cO_\Delta) \ra \Hom_{X\times
  X}(\cO_\Delta, S_\Delta) \]
is the isomorphism given by 
\[ \mu \mapsto \tau(\mu) = \mu\otimes \pi_2^* S_X, \]
and we considered a similar isomorphism
\[ \bar{\tau}:\Hom_{X\times X}(S_\Delta^{-1}, \cO_\Delta) \ra \Hom_{X\times
  X}(\cO_\Delta, S_\Delta) \]
given by
\[ \mu \mapsto \bar{\tau}(\mu) = \mu\otimes \pi_1^* S_X, \]
which corresponds to choosing left adjoints in the definition of the
Mukai product, instead of right adjoints, as we did in
Definition~\ref{def:mukprod}.  We extend this notation to simply mean
that $\tau$ is the operation of tensoring with $\pi_2^* S_X$, and
$\bar{\tau}$ is the similar operation that corresponds to $\pi_1^*
S_X$.

Let 
\begin{align*}
\alpha = \tau(v) & :\cO_{\Delta_Y}\ra S_{\Delta_Y}, \\
\beta = \bar{\tau}(w) & : \cO_{\Delta_X} \ra S_{\Delta_X},
\end{align*}
and consider the morphisms
\[
\begin{diagram}[height=2em,width=2em,labelstyle=\scriptstyle]
\mu_1 : \cF & \rTo^{\cF\circ \beta} & \cF \circ S_{\Delta_X} &
  \rTo^{\alpha\circ \cF\circ S_{\Delta_X}} & S_{\Delta_Y} \circ
  \cF\circ S_{\Delta_X}, \\
\mu_2 : \cF & \rTo^{\alpha \circ \cF} & S_{\Delta_Y} \circ \cF &
  \rTo^{S_{\Delta_Y} \circ \cF\circ \beta} & S_{\Delta_Y} \circ
  \cF\circ S_{\Delta_X}.
\end{diagram}
\]
By the commutativity of the trace (Lemma~\ref{lem:shuffleok}) it follows that 
\[ \Tr_{X\times Y}(\mu_1) = \Tr_{X\times Y}(\mu_2). \]
Now consider the commutative diagrams
%\[ 
%\begin{diagram}[height=2em,width=2em,labelstyle=\scriptstyle]
%\cO_{\Delta_X} & \rTo^{\bar{\eta}} & \cH\circ \cF & & & \cO_{\Delta_Y}
%& \rTo^{\eta} & \cF\circ \cG \\
%\dTo^\beta & & \dTo_{\cH\circ \cF\circ \beta} & & & \dTo^\alpha & &
%\dTo_{\alpha\circ \cF\circ \cG} \\
%S_{\Delta_X} & \rTo^{\bar{\eta}\circ S_{\Delta_X}} & \cH\circ\cF\circ
%  S_{\Delta_X} & & & S_{\Delta_Y} & \rTo^{S_{\Delta_Y}\circ \eta} &
%  S_{\Delta_Y}\circ \cF\circ \cG \\
%& & \dTo_{\cH\circ \alpha\circ \cF\circ S_{\Delta_X}} & & & & &
%  \dTo_{S_{\Delta_Y}\circ \cF\circ\beta\circ \cG} \\
%& & \cH\circ S_{\Delta_Y}\circ
%\cF S_{\Delta_X} & & & & &
%S_{\Delta_Y}\circ \cF\circ S_{\Delta_X} \circ \cG \\
%\dTo^{(\Phi^*\alpha)\circ S_{\Delta_X}} & & \dEqual & & & \dTo^{S_{\Delta_Y}\circ (\Psi^*\beta)} & & \dEqual \\
%& & S_{\Delta_X}\circ\cG\circ\cF\circ S_{\Delta_X} & & & & &
%S_{\Delta_Y}\circ\cF\circ\cH\circ S_{\Delta_Y} \\
%& & \dTo_{S_{\Delta_X}\circ \epsilon \circ S_{\Delta_X}} & & & & &
%\dTo_{S_{\Delta_Y} \circ \bar{\epsilon} \circ S_{\Delta_Y}} \\
%S_{\Delta_X}\circ S_{\Delta_X} & \rEqual & S_{\Delta_X}\circ
%S_{\Delta_X} & & &  S_{\Delta_Y}\circ S_{\Delta_Y}& \rEqual &
%S_{\Delta_Y}\circ S_{\Delta_Y}. 
%\end{diagram}
%\]
\[ 
\begin{diagram}[height=2em,width=2em,labelstyle=\scriptstyle]
1_{X} & \rTo^{\bar{\eta}} & \cH \cF & & & 1_{Y}
& \rTo^{\eta} & \cF \cG \\
\dTo^\beta & & \dTo_{\cH \cF \beta} & & & \dTo^\alpha & &
\dTo_{\alpha \cF \cG} \\
S_{X} & \rTo^{\bar{\eta} S_{X}} & \cH\cF
  S_{X} & & & S_{Y} & \rTo^{S_{Y} \eta} &
  S_{Y} \cF \cG \\
& & \dTo_{\cH \alpha \cF S_{X}} & & & & &
  \dTo_{S_{Y} \cF\beta \cG} \\
& & \cH S_{Y}
\cF S_{X} & & & & &
S_{Y} \cF S_{X}  \cG \\
\dTo^{(\Phi^\dagger\alpha) S_{X}} & & \dEqual & & & \dTo^{S_{Y} (\Psi^\dagger\beta)} & & \dEqual \\
& & S_{X}\cG\cF S_{X} & & & & &
S_{Y}\cF\cH S_{Y} \\
& & \dTo_{S_{X} \epsilon  S_{X}} & & & & &
\dTo_{S_{Y}  \bar{\epsilon}  S_{Y}} \\
S_{X} S_{X} & \rEqual & S_{X}
S_{X} & & &  S_{Y} S_{Y}& \rEqual &
S_{Y} S_{Y},
\end{diagram}
\]
where we have omitted the $\circ$ signs, and we wrote $S_X$ for
$S_{\Delta_X}$ and $1_X$ for $\cO_{\Delta_X}$.  

Reading around the diagrams and using Proposition~\ref{prop:xytoxx} we
see that
\begin{align*}
 \Tr_{X\times X}((\Phi^\dagger\alpha)S_X\circ \beta) & = \Tr_{X\times
  X}(S_X\epsilon S_X \circ \cH\mu_1\circ \bar{\eta}) \\
& = \Tr_{X\times Y}(\mu_1) = \Tr_{X\times Y}(\mu_2) \\
& = \Tr_{Y\times Y}(S_Y\bar{\epsilon}S_Y
  \circ \mu_2\cG \circ \eta) \\
& = \Tr_{Y\times Y}(S_Y(\Psi^\dagger\beta)\circ \alpha).
\end{align*}
Reverting to the $\tau$, $\bar{\tau}$ notation (where multiplication
by $S_X$ or $S_Y$ on the left corresponds to $\tau$, and on the right
to $\bar{\tau}$), we conclude that
\[ \Tr_{X\times X}(\bar{\tau}\Phi^\dagger\tau v \circ \bar{\tau} w) =
\Tr_{Y\times Y}(\tau\Psi^\dagger\bar{\tau}w \circ \tau v), \]
or, since $\tau$, $\bar{\tau}$ are simply multiplication by a line bundle 
\begin{align*}
\Tr_{X\times X}(\Phi^\dagger \tau v \circ w) & = \Tr_{Y\times
  Y}(\Psi^\dagger\bar{\tau} w \circ v) \\ 
& = \Tr_{Y\times Y}(\bar{\tau}w \circ  \Psi_* v) \\
& = \Tr_{Y\times Y}(\tau\Psi_* v \circ w),
\end{align*}
where the second equality is the definition of $\Psi_*$, and the third
one follows from the fact that $\tau\bar{\tau} = S_{Y\times Y}$ and
Lemma~\ref{lem:shuffleok}.

Since $w$ was arbitrary and the pairings are non-degenerate, we
conclude that $\Phi^\dagger\tau = \tau\Psi_*$, and this completes the proof.
\end{proof}

\remark
Although we have 
\begin{align*}
\tau(\cO_{\Delta_X}) & = \bar{\tau}(\cO_{\Delta_X}) = S_{\Delta_X},
\intertext{and}
\tau(S_{\Delta_X}) & = \bar{\tau}(S_{\Delta_X}) = S_{X\times X}
\cO_{\Delta_X},
\end{align*}
the two isomorphisms
\[ \tau,\ \bar{\tau}  : \Hom_{X\times X}(\cO_{\Delta_X}, S_{\Delta_X})
\ra \Hom_{X\times X}(S_{\Delta_X}, S_{X\times X} \cO_{\Delta_X}) \]
are {\em different}.  This can essentially be seen by looking at Chern
characters.  This is the reason we need to be careful about the
distinction between $\tau$ and $\bar{\tau}$.

\begin{corollary}
\label{cor:isoisom}
If $\Phi:\D(X)\ra \D(Y)$ is an equivalence, then $\Phi_*:HH_*(X)\ra
HH_*(Y)$ is an isometry with respect to the generalized Mukai product.
\end{corollary}

\begin{proof}
We have
\[ \langle\, \Phi_* x, \Phi_* y\,\rangle = \langle\, \Psi_*\Phi_* x,
y\,\rangle = \langle\, (\Psi\circ \Phi)_* x, y \,\rangle = \langle\,
x,y\,\rangle, \]
where the last equality follows from the fact that if $\Phi$ is an
equivalence, then its left adjoint $\Psi$ is an inverse to it.
\end{proof}

\subsection{}
The Hirzebruch-Riemann-Roch theorem is a consequence of the other
properties:
\begin{theorem}
\label{thm:hrr}
For $\cE, \cF\in \D(X)$ we have 
\[ \MP{\ch(\cE)}{\ch(\cF)} = \chi(\cE, \cF), \]
where $\chi(\scdot,\scdot)$ is the Euler pairing on $K_0(X)$, 
\[ \chi(\cE, \cF) = \sum_i (-1)^i \dim \Ext^i_X(\cE, \cF). \]
\end{theorem}

\begin{proof}
Let $p:X\ra \pt$ be the structure morphism of $X$, and observe that
$\cO_X = p^* \cO_\pt$.  The functor $p^*$ is left adjoint to $p_*$,
and if $\Phi$ is the functor $\cE\otimes\,-$, then its right adjoint
$\Psi$ is given by $\cE^\chk\otimes\, -$.

Using the properties of the Mukai product and Chern character we get
\begin{align*}
\MP{\ch(\cE)}{\ch(\cF)} & = \MP{\ch(\cE\otimes\cO_X)}{\ch(\cF)} \\
& = \MP{\ch(\Phi\cO_X)}{\ch(\cF)} \\
& = \MP{\Phi_*\ch(\cO_X)}{\ch(\cF)} \\
& = \MP{\ch(\cO_X)}{\Psi_*\ch(\cF)} \\
& = \MP{\ch(\cO_X)}{\ch(\Psi\cF)} \\
& = \MP{\ch(\cO_X)}{\ch(\cE^\chk \otimes \cF)} \\
& = \MP{\ch(p^* \cO_\pt)}{\ch(\cE^\chk \otimes \cF)} \\
& = \MP{(p^*)_*\ch(\cO_\pt)}{\ch(\cE^\chk \otimes \cF)} \\
& = \MP{\ch(\cO_\pt)}{(p_*)_*\ch(\cE^\chk \otimes \cF)} \\
& = \MP{\ch(\cO_\pt)}{\ch(p_*(\cE^\chk \otimes \cF))}\\
& = \MP{\ch(\cO_\pt)}{\ch(\R\Hom_X(\cE, \cF))}.
\end{align*}
Since $\ch$ is a map on K-theory, $K_0(\pt) \iso \Z$, and the Mukai
product is additive, we see that
\[   \MP{\ch(\cO_\pt)}{\ch(\R\Hom_X(\cE, \cF))} =
\chi(\R\Hom_X(\cE, \cF))\scdot \MP{\ch(\cO_\pt)}{\ch(\cO_\pt)} =
\chi(\cE, \cF), \]
as it is a trivial computation to check that 
\[ \MP{\ch(\cO_\pt)}{\ch(\cO_\pt)}_\pt = 1. \]
\end{proof}

\remark
The same proof works in the case of a global quotient orbifold, with
some minor corrections.  The map $p: X \ra \pt$ needs to be replaced
by $p:X\ra [\pt/G] = BG$.  Then we consider $\cO_\pt$ with the trivial
$G$-representation, and $p^* \cO_\pt$ is what we use for $\cO_X$.  The
K-theory of $BG$ is not $\Z$, but we are only interested in the part
that is generated by $\cO_\pt$ (with the trivial representation),
since we are only interested in $G$-equivariant $\Hom$ groups in the
definition of $\chi(\cE, \cF)$ (see~\cite[Section 4]{BKR} for
details).

\remark
If we define $\Td(X) = \tau\ch(\cO_X)$, we have the following
alternative version of the above theorem
\[ \chi(\cF) = \Tr_{X\times X}(\Td(X) \circ \ch(\cF)), \]
more reminiscent of the classical version of Hirzebruch-Riemann-Roch.

\subsection{}
We conclude with a mention of the following result, inspired by the
Cardy condition in physics.  We omit the proof, as we shall not use it
in the sequel, and it is mainly an exercise in applying several times
the basic construction (Proposition~\ref{prop:fundprop}).  The
interested reader can easily supply the details.

\begin{theorem}[Cardy condition]
\label{thm:cardy}
Let $\cE, \cF$ be objects in $\D(X)$, and let $e\in \Hom_X(\cE, \cE)$
and $f\in \Hom_X(\cF, \cF)$.  Consider the operator
\[ {}_f m_e :\Hom^i_X(\cE, \cF) \ra \Hom^i_X(\cE, \cF) \]
given by composition by $f$ on the left and by $e$ on the right.  Then
we have
\[ \MP{\iota^\cE\!(e)}{\iota^\cF\!(f)} = \sTr {}_f m_e, \]
where $\iota^\cE$, $\iota^\cF$ are the maps defined
in~(\ref{subsec:iotae}), and $\sTr$ denotes the alternating sum of the
traces of the action of ${}_f m_e$ on $\Hom_X^i(\cE, \cF)$.
\end{theorem}

Observe that the Hirzebruch-Riemann-Roch formula is a direct
consequence of the Cardy condition, with $e=\id_\cE$, $f=\id_\cF$.

\section{Derived equivalence invariance}
\label{sec:invderequiv}

This section is devoted to a discussion of the invariance of the
Hochschild structure under derived equivalences.  This is the primary
reason for our decision to use it instead of the harmonic structure
given by cohomology of vector fields and/or forms discussed
in~\cite{CalHH2}.  We provide proofs of our statements in the
``$\Ext$'' interpretation; it is obvious that the proofs in the
categorical interpretation would be significantly shorter, perhaps
trivial.

\subsection{}
We aim to prove the following result:
\begin{theorem}
\label{thm:dcatinv}
Let $X$ and $Y$ be spaces whose derived categories are equivalent via
a Fourier-Mukai transform (i.e., the equivalence is given by an
integral transform).  Then there exists a natural isomorphism of
Hochschild structures
\[ (HH^*(X), HH_*(X)) \iso (HH^*(Y), HH_*(Y)). \]
\end{theorem}

\subsection{}
There are three statements implicit in the above theorem, which we'll
discuss in turn:
\begin{itemize}
\item[a.] $HH^*(X) \iso HH^*(Y)$ as graded rings;
\item[b.] $HH_*(X) \iso HH_*(Y)$ as graded modules over the cohomology
  rings;
\item[c.] the isomorphism $HH_*(X)\iso HH_*(Y)$ is an isometry with
  respect to the generalized Mukai product.
\end{itemize}

\begin{proposition}
\label{prop:eqvdeltaS}
Under the hypothesis of Theorem~\ref{thm:dcatinv} there is an
equivalence of derived categories 
\[ \D(X\times X) \iso \D(Y\times Y) \]
which maps $\cO_{\Delta_X}$ to $\cO_{\Delta_Y}$ and $S_{\Delta_X}$ to
$S_{\Delta_Y}$.
\end{proposition}

Again, this statement would be trivial in the ``natural
transformations'' context: $\cO_{\Delta_X}$ and $\cO_{\Delta_Y}$
correspond to the identity natural transformations, and $S_{\Delta_X}$
and $S_{\Delta_Y}$ correspond to the Serre functors (intrinsic to any
triangulated category that possesses one).

\begin{proof}
(Independently also proven in~\cite{GolLunOrl},~\cite{OrlAb}.)  Begin
with the observation that if $\cE$ is an object in $\D(X\times Y)$
that induces a Fourier-Mukai transform $F=\FMXY^\cE:\D(X)\ra \D(Y)$,
then $F^\chk = \FMXY^{\cE^\chk}$ is also an equivalence.  (We have
{\em not} interchanged $X$ and $Y$, as expected.)  Indeed,
\[ F^\chk\cO_x = (F\cO_x)^\chk \]
for $x\in X$, and thus
\begin{align*}
\Hom_Y(F^\chk\cO_x, F^\chk\cO_{x'}) & =
\Hom_Y((F\cO_x)^\chk, (F\cO_{x'})^\chk) \\
& = \Hom_Y(F\cO_{x'}, F\cO_x) 
\end{align*}
for $x,x'\in X$.  It follows that the orthogonality condition
of~\cite[Theorem 5.1]{BriFM} is satisfied by $F^\chk$ if it is already
satisfied by $F$.  Furthermore, if $F$ is an equivalence
\[ F^\chk\cO_x = (F\cO_x)^\chk \iso (F\cO_x \otimes \omega_Y)^\chk =
(F\cO_x)^\chk \otimes \omega_Y^{-1} = F^\chk\cO_x \otimes
\omega_Y^{-1} \]
by~\cite[Theorem 5.4]{BriFM}, and therefore
\[ F^\chk\cO_x\otimes \omega_Y \iso F^\chk\cO_x, \]
so we conclude by the same theorem that $F^\chk$ is an equivalence.

Now consider the objects $\cE^*=\cE^\chk\otimes \pi_Y^* S_Y$ and $\cE$
on $X\times Y$.  Both induce equivalences $\D(X)\ra\D(Y)$, so
by~\cite[Corollary 1.8]{OrlAb} we conclude that if we let 
\[ \cH = \cE^*\boxtimes \cE = \pi_{12}^* \cE^* \otimes \pi_{34}^* \cE
\]
then $\cH$ induces an equivalence of derived categories
$\Phi^\cH:\D(X\times X) \ra \D(Y\times Y)$.  (Here, $\pi_{ij}:X\times
Y\times X\times Y\ra X\times Y$ are the projections onto the
corresponding factors.)

We claim that $\Phi^{\cH}(\cO_{\Delta_X}) = \cO_{\Delta_Y}$.  Indeed,
we have
\[ \Phi^{\cH}(\cO_{\Delta_X}) = \pi_{24,*}(\pi_{13}^*(\cO_{\Delta_X})
\otimes \cH), \]
and 
\[ \pi_{13}^* \cO_{\Delta_X} = \Delta_{13,*}\cO_{Y\times X\times Y}, \]
where $\Delta_{13}$ is the morphism
\[ \Delta_{13}:Y\times X\times Y \mapsto X\times Y\times X\times Y \]
which maps
\[ (y_1,x,y_2) \mapsto (x,y_1,x,y_2). \]
Therefore, by the projection formula we have 
\[ \pi_{13}^*(\cO_{\Delta_X}) \otimes \cH = \Delta_{13,*}
\Delta_{13}^* \cH, \]
and by the commutative diagram
\[ 
\begin{diagram}[height=2em,width=4em,labelstyle=\scriptstyle,silent]
Y\times X\times Y & \rTo^{\ \ \ \Delta_{13}\ \ \ } & X\times Y\times X\times Y \\
\dTo^{p_{13}} & \ldTo_{\pi_{24}} & \\
Y\times Y,
\end{diagram}
\]
we conclude that
\begin{align*}
\Phi^{\cH}(\cO_{\Delta_X}) & = \pi_{24,*}(\pi_{13}^*(\cO_{\Delta_X})
\otimes \cH) \\
& = \pi_{24,*}\Delta_{13,*} \Delta_{13}^* \cH  = p_{13,*}
\Delta_{13}^* \cH \\
& = p_{13,*}( \cE^* \boxtimes' \cE),
\end{align*}
where we have denoted by $p_{ij}$ the projections from $Y\times
X\times Y$ onto the corresponding factors, and 
\[ \cE^*\boxtimes'\cE = p_{12}^* \cE^* \otimes p_{23}^*\cE. \]

But with notation as in Section~\ref{sec:prel},
\[ p_{13,*}(\cE^*\boxtimes'\cE) = \cE \circ \cE^*, \]
and since $\FMXY^\cE$ and $\FMYX^{\cE^*}$ were inverse to one another,
\[ \cE\circ\cE^* \iso \cO_{\Delta_Y} \]
by the comments in~(\ref{subsec:composite}).  We conclude that
\[ \Phi^{\cH}(\cO_{\Delta_X}) = \cO_{\Delta_Y}. \]

The computation of $\Phi^{\cH}(S_{\Delta_X})$ is entirely similar, and
we shall omit the details.  We shall only mention that what one
obtains is that 
\[ \Phi^{\cH}(S_{\Delta_X}) = \cE \circ S_{\Delta_X} \circ \cE^*, \]
and in terms of functors that corresponds to 
\[ F\circ S_X\circ F^{-1}, \]
where $F=\Phi^\cE$.  But since the Serre functor is
intrinsic~\cite{BonKapSerre}, this functor must be isomorphic to
$S_Y$, and hence it must be given by $S_{\Delta_Y}$ (an equivalence
of derived categories is induced by an object on the product which is
unique up to isomorphism~\cite[Theorem 2.2]{Orl}).
\end{proof}

\begin{corollary}
A Fourier-Mukai transform $\D(X)\iso \D(Y)$ induces an isomorphism of
graded rings $HH^*(X) \iso HH^*(Y)$, as well as an isomorphism
$HH_*(X) \iso HH_*(Y)$ of graded modules over the corresponding
cohomology rings.
\end{corollary}

\begin{proof}
Follows immediately from Proposition~\ref{prop:eqvdeltaS}.  It is
useful to point out that this isomorphism is independent of the choice
of isomorphism $\Phi^\cH(\cO_{\Delta_X})\iso \cO_{\Delta_Y}$: indeed,
different choices are conjugate by an element of $\Hom_{Y\times
  Y}(\cO_{\Delta_Y}, \cO_{\Delta_Y})$, and these elements are central
in $HH^*(Y)$.  A similar argument works for $HH_*$.

In the affine case this theorem has been proven for cohomology even
for non-commutative rings by Happel~\cite{Hap} and
Rickard~\cite{Ric}.  (Rickard even removed the requirement that the
equivalence be given by a Fourier-Mukai transform.)
\end{proof}

This completes parts a.\ and b.\ of Theorem~\ref{thm:dcatinv}.  Part
c.\ is just Corollary~\ref{cor:isoisom}.

\appendix

\section{A categorical approach via topology and TQFT's}
\label{app:categ1}

In this appendix we present two categorical/topological/TQFT
approaches to the Hochschild structure.  Some of the ideas will not be
completely rigorous, however the intuition behind these approaches is
often extremely valuable in understanding the proofs in this paper.

\subsection{}
Consider the weak 2-category $\cD$ (the 2-category of ``all derived
categories''), defined as follows:
\begin{itemize}
\item[--] objects of $\cD$ are smooth, projective schemes (or compact
  orbifolds, or twisted spaces, etc.); we'll call the objects of $\cD$
  {\em spaces};
\item[--] if $X$ and $Y$ are spaces, $\Hom(X, Y) = \Ob
  \D(X\times Y)$; 
\item[--] if $\cE, \cF$ are elements of $\Hom(X, Y)$, then $\THom(\cE,
  \cF) = \Hom_{\D(X\times Y)}(\cE, \cF)$.
\end{itemize}
The composition of 1-morphisms is given by convolution of kernels, as
in~(\ref{subsec:composite}).  The horizontal and vertical composition
of 2-morphisms is defined in the obvious way.  The (weak) identity
1-morphism $X\ra X$ is given by $\cO_{\Delta_X}\in \Hom(X, X) =
\D(X\times X)$.  (The reader unfamiliar with 2-categories is referred
to~\cite{Bae}.)

Observe that $\cD$ has a richer structure than just that of a 2-category:
if $X$ and $Y$ are spaces, $\Hom(X, Y)$ has a natural structure of a
triangulated category which admits a Serre functor.

\subsection{}
It is useful to consider the 2-functor $\Phi:\cD \ra \gCat$
between the 2-category $\cD$ and the 2-category $\gCat$ of all
categories (where 1-morphisms are functors and 2-morphisms are natural
transformations).  The functor $\Phi$ is defined by setting
\begin{enumerate}
\item[--] $\Phi(X) = \D(X)$ for any space $X$;\vspace{0.5mm}
\item[--] $\Phi(\cE) = \FMXY^\cE:\D(X)\ra\D(Y)$ for any object
  $\cE\in\D(X\times Y) = \Hom_\cD(X, Y)$;\vspace{0.5mm}
\item[--] $\Phi(f) = \FMXY^f:\FMXY^\cE \Rightarrow \FMXY^\cF$ for any
  morphism $f:\cE\ra \cF$ between $\cE,\cF\in\D(X\times Y)$.\vspace{0.5mm}
\end{enumerate}
Note again that the natural properties associated with the
triangulated structure are preserved ($\Phi(\cE)$ is exact, $\Phi(f)$
commutes with translations, etc.)  Example~\ref{example:counter} shows
that $\Phi$ can not be fully faithful, even if we impose extra
conditions on what functors or natural transformations we allow in
$\gCat$.  In spite of this, we shall think of $\Phi$ as being fully
faithful, and thus think of morphisms of $\cD$ as functors,
2-morphisms of $\cD$ as natural transformations.  For more on this
problem and a possible way around it see Appendix~\ref{app:categ2}.

\subsection{}
Our purpose is to find invariants of spaces that are, in a sense to be
made precise later, functorial with respect to $\cD$.  As a first
example of such an invariant, consider associating to a space $X$ the
group $K_0(X\times X)$.  The functoriality of $K_0$ is expressed by
the fact that for a morphism $X\ra Y$ in $\cD$ given by an object
$\cE\in\D(X\times Y)$, there is a natural map $K_0(X\times X)\ra
K_0(Y\times Y)$ defined by
\[ \cF\in\D(X\times X) \mapsto \cE \circ \cF \circ \cE^* \]
where 
\[ \cE^* = \cE^\chk \otimes \pi_Y^* S_Y, \]
is the object in $\D(Y\times X)$ whose associated functor $\D(X) \ra
\D(Y)$ is left adjoint of $\FMXY^\cE$.  It is easy to see that this
association of morphisms in $\cD$ and maps between $K_0$ groups is
functorial.

\subsection{}
A useful analogy is obtained by considering the 2-category $\cT$
associated to a topological space $T$, where objects of $\cT$ are
points of $T$, morphisms between two points $x,y\in T$ are given by
paths in $T$ from $x$ to $y$, and 2-morphisms are given by homotopy
equivalence classes of homotopies between paths.  Observe that in
$\cT$ all morphisms are isomorphisms.

The first homotopy invariant of $T$ that one studies is the fundamental
group $\pi_1(T, x)$ of homotopy classes of loops based at $x\in T$.
Categorically, this can be thought of as the set of morphisms $x\ra x$
in $\cT$ (loops at $x$), modulo the equivalence relation induced
by 2-isomorphisms (i.e., homotopies).  Given a morphism $x\ra y$
(a path from $x$ to $y$), one obtains a natural map $\pi_1(T,x)
\ra \pi_1(T,y)$ by conjugating a loop at $x$ with the given path $x\ra
y$.  This map only depends on the homotopy class of the path $x\ra y$,
and the association of maps to paths is functorial.

It is now obvious that the same procedure that was used to construct
the fundamental group of $T$ has also been used to construct
$K_0(X\times X)$ in $\cD$.  (In fact, the analogy is imperfect: in
constructing $K_0(X\times X)$ we have taken the quotient of
$\Ob \D(X\times X)$ by a far larger equivalence relation: while in
defining $\pi_1$ we only identified 1-morphisms that were 2-isomorphic,
in constructing $K_0(X\times X)$ we have declared that 1-morphisms
that form a triangle should sum to zero.  But in the context of the
existence of the triangulated structure on $\Hom_\cD(X, Y)$ it makes
sense to consider this coarser equivalence relation to get a more
finite invariant.)  The analogy also extends to the association of
morphisms between $K_0$-groups to 1-morphisms in the underlying
category.

The above discussion also shows one of the most important weaknesses
of this analogy: while in $\cT$ every morphism $x\ra y$ is an
isomorphism, and hence it makes sense to talk about its inverse when
conjugating a loop at $x$ to get a loop at $y$, morphisms in $\cD$ are
not invertible in general.  Thus we had to settle for the weaker
concept of left adjoint; but there was no particular reason to choose
left over right: right adjoints would have worked equally well.

\subsection{}
The procedure just described should be thought of as producing a
functor 
\[ K_0: \OcD \ra \gGps \]
from the decategorification $\OcD$ of $\cD$ (the 1-category that is
obtained from $\cD$ by forgetting 2-morphisms and setting 1-morphisms
that are 2-isomorphic in $\cD$ to be equal in $\OcD$) to the category
$\gGps$ of groups.  There is no reason, however, to stop at the
$\pi_1$ level in our analogy with topological spaces.  Indeed, note
that $\pi_2(T, x)$, as the space of homotopy classes of maps $(S^2,
\pt) \ra (T, x)$, can be thought of as the space of homotopy classes
of homotopies from the constant path $\Id_x$ at $x$ to itself.  In
other words
\[ \pi_2(T, x) = \THom_\cT(\Id_x, \Id_x). \]
Just as in the case of $\pi_1$, a path $x\ra y$ induces in a
functorial way a map $\pi_2(T,x) \ra \pi_2(T, y)$.

\subsection{}
\label{subsec:defhh}
By analogy, in the category $\cD$ we associate to a
space $X$ its Hochschild cohomology
\[ HH^*(X) = \THom_\cD(\Id_X, \Id_X) = \Hom_{\D(X\times
  X)}(\cO_\Delta, \cO_\Delta). \] 
However, when we try to mimic the construction of the map $\pi_2(T, x)
\ra \pi_2(T, y)$ associated to a map $x\ra y$ to obtain a map
$HH^*(X)\ra HH^*(Y)$ associated to a map (in $\cD$) $X\ra Y$ we hit a
major difficulty: in the topological setting, the construction relies
on the fact that if $f:x\ra y$ is a path in $T$, then $f^{-1} \circ f$
is homotopic to the constant path $\Id_x$, {\em and this homotopy can be
read as either a 2-morphism $1\Rightarrow f^{-1} \circ f$ or
$f^{-1}\circ f \Rightarrow 1$}.  In $\cD$ this is no longer the case:
we have to replace a path $x\ra y$ by an object $\cE\in\D(X\times Y)$,
and there is no reason why the functor $\FMXY^\cE$ should have an
inverse.  It will always have left and right adjoints, but these need
not be isomorphic in general.  In order to compensate for this
discrepancy we need to modify the definition of $HH^*$ to
\[ HH_*(X) = \THom_\cD(S_X^{-1}, \Id_X) = \Hom_{\D(X\times
  X)}(\Delta_! \cO_X, \cO_\Delta). \] 
As it turned out in Section~\ref{sec:functoriality}, with this definition it
is possible to mimic the construction that we did above for $\pi_2$,
and to get a functor
\[ HH_*:\OcD \ra \gVect \]
from the decategorification of $\cD$ to the category of graded vector
spaces.

One can also see that when dealing with an isomorphism $X\ra Y$ in
$\cD$ (i.e., and equivalence of categories) the above technique can be
used to construct an isomorphism $HH^*(X)\iso HH^*(Y)$ (we use the
fact that the left and right adjoints of an equivalence are
isomorphic).  This is precisely what we did in
Section~\ref{sec:invderequiv}.

\subsection{}
There is a second approach to Hochschild homology via topology, which
in a sense is orthogonal to the above one.  In this approach we
consider the 2-category $\cC$ of 3-cobordisms with corners defined as
follows (for details see~\cite{Lau}):
\begin{itemize}
\item[--] objects of $\cC$ are smooth, oriented 1-manifolds (i.e.,
  disjoint unions of circles);
\item[--] 1-morphisms of $\cC$ are cobordisms between the objects of
  $\cC$ (smooth oriented surfaces with boundary);
\item[--] 2-morphisms of $\cC$ are cobordisms between 1-morphisms
  (3-manifolds with corners).
\end{itemize}

We are interested in studying 2-functors from $\cC\ra \cD$.  Such
functors, when they satisfy certain other properties (the list of
these properties varies in the literature), are also known as extended
TQFT's (or 1+1+1 TQFT's).  One of the more common requirements on
TQFT's is that they should respect the monoidal structure on objects,
given in $\cC$ by disjoint union of 1-manifolds and in $\cD$ by
product of spaces, and we shall search for functors with this
property.  When a TQFT $F:\OcC\ra\OcD$ is only defined between the
decategorifications $\OcC$, $\OcD$ of $\cC$ and $\cD$, it is known as
a 1+1 TQFT.

\subsection{}
Building on work of Roberts, Sawon and Willerton (unpublished) we
discuss a possible approach to constructing such a functor (their
construction is essentially the same as ours, but with a different
target category).  Due to technical difficulties we shall actually
only discuss their construction of a 1-functor $\cX:\OcC\ra\OcD$; the
main point of our discussion is to argue that despite the fact that we
can prove that this functor can not be extended to a 2-functor $\cC\ra
\cD$, the functoriality property of Hochschild homology should be
viewed as strong evidence that there exists a modification of the
category $\cC$ for which such a lifting exists.  Such a theory should
yield interesting invariants of 3-manifolds.

\subsection{}
Let us begin first with the construction of a 1-functor $\cX:\OcC\ra
\OcD$, which depends on the choice of a space $X$.  Associate to $S^1$
the space $X$, and since we want $\cX$ to respect the monoidal
structure on objects,
\[ \cX(S^1 \coprod \cdots \coprod S^1) = X\times \cdots \times X, \]
so in particular $\cX(\emptyset) = \pt$.  To a pair of pants associate
$\cO_\Delta\in\D(X\times X\times X)$, where $\Delta$ is the small
diagonal 
\[ \{(x,x,x)\in X\times X\times X~|~x\in X\} \]
in the product $X\times X\times X$ (we think of objects in
$\D(X\times X\times X)$ as morphisms $X\ra X\times X$ in $\cD$).  As
functors, this corresponds either to 
\[ \Delta_*:\D(X)\ra \D(X\times X) \]
or to 
\[ \Delta^*:\D(X\times X) \ra \D(X), \]
depending on whether we want to think of the pair of pants as a
cobordism from $S^1$ to $S^1\coprod S^1$, or the other way around.
We associate to the disk (thought of as a cobordism $\emptyset\ra
S^1$) the object $\cO_X\in \D(\pt \times X)$.  As a functor it
corresponds either to $p^*$ or to $p_*$ depending on whether we view
it as a map $\emptyset \ra S^1$ or $S^1\ra \emptyset$ ($p:X\ra \pt$ is
the structure map of $X$).

Knowing this information is enough to determine the value of $\cX$ on
any 1-cobordism (any oriented surface with boundary can be decomposed
into a finite number of pairs of pants and caps).
\begin{figure}
\begin{center}
\includegraphics*{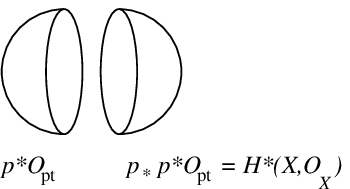}
\caption{Computation of $\cX(S^2)$}
\label{fig:fig1}
\end{center}
\end{figure}
Figures~\ref{fig:fig1} and~\ref{fig:fig2} show how to compute
$\cX(S^2)$ and $\cX(T^2)$ (we think of a closed surface as a cobordism
from the empty manifold to itself, and thus $\cX(\mbox{closed
surface})$ is a map from a point to itself in $\cD$, which is just a
graded vector space).
\begin{figure}
\begin{center}
\includegraphics*{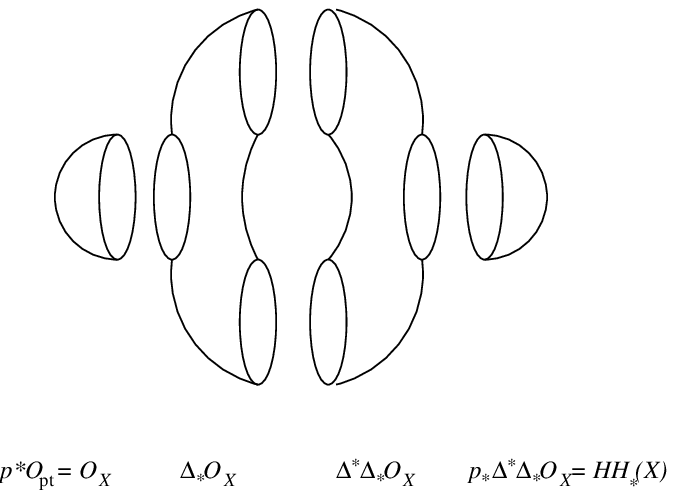}
\caption{Computation of $\cX(T^2)$}
\label{fig:fig2}
\end{center}
\end{figure}
The results are
\begin{align*}
\cX(S^2) & = p_* \cO_X = H^*(X, \cO_X) \\
\cX(T^2) & = p_* \Delta^*\Delta_* \cO_X = \R\Hom_X(\cO_X,
\Delta^*\Delta_*\cO_X) \\
& = \R\Hom_{X\times X}(\Delta_! \cO_X, \Delta_* \cO_X) = HH_*(X).
\end{align*}
Thus we conclude that this 1+1 TQFT with target $X$ should associate
to $T^2$ the Hochschild homology $HH_*(X)$ of $X$.

\subsection{}
Unfortunately it is known that $\cX$ can not be lifted to a 2-functor
$\cX:\cC\ra\cD$.  The obstruction to doing this has already appeared
in the work of Khovanov~\cite{Kho}: his observation is that if one has
a 2-functor $F:\cC\ra \gCat$, then if $S$ is any 2-dimensional
cobordism (1-morphism in $\cC$) between the 1-manifolds $C_1$ and
$C_2$, then $F(S)$, which is a functor between the categories $F(C_1)$
and $F(C_2)$ must have a biadjoint (a functor which is both a
left and a right adjoint to $F(S)$).  This biadjoint is given by
$F(S')$, where $S':C_2\ra C_1$ is the reverse cobordism, given by the
same manifold $S$.

In particular this shows that the functor $\cX$ that we constructed
earlier has no hope of lifting to a 2-functor $\cC\ra \cD$ unless the
target space is a point: the composition $\Phi\circ \cX$ would be a
2-functor $\cC\ra \gCat$, which associates to a pair of pants (thought
of as a morphism $S^1\ra S^1\coprod S^1$) the functor $\Delta_*:\D(X)
\ra \D(X\times X)$, and to the reverse pair of pants 
$\Delta^*:\D(X\times X) \ra \D(X)$.  Unfortunately, while $\Delta^*$
is a left adjoint to $\Delta_*$, it is not a right adjoint unless
$X=\pt$.

\subsection{}
Let us, however, assume that an extension of $\cX$ to a 2-functor
could be found, for every target space $X$.  Observe that once the
image of $\cX$ on $S^1$ is fixed, the entire functor $\cX$ is
essentially fixed.  We can expect that a similar statement would hold
for natural transformations between such functors: given functors
$\cX,\cY:\cC\ra\cD$ that correspond to spaces $X, Y$, respectively, a
natural transformation (with certain properties yet to be fixed)
should be completely determined by its value on $S^1$.  But its value
on $S^1$ is just a map $\cX(S^1)\ra \cY(S^1)$ in $\cD$, i.e., an
object in $\D(X\times Y)$.  Thus given spaces $X$ and $Y$, we should
get associated 2-functors $\cX$ and $\cY$, and given an object
$\cE\in\D(X\times Y)$, we should get an associated natural
transformation between them.

Observe that a natural transformation $\eta:\cX\ra\cY$ of 2-functors
between $\cC$ and $\cD$ is a collection of 1-morphisms $\eta(o):\cX(o)
\ra \cY(o)$ associated to objects $o\in\Ob\cC$, as well as a
collection of 2-morphisms $\eta(o\ra o')$ associated to 1-morphisms in
$\cC$.  The 2-morphism $\eta(o\ra o')$ is depicted in the following
diagram
\[
\begin{diagram}[height=1.25em,width=4em,labelstyle=\scriptstyle]
& & \cX(o') & & \\
& \ruTo^{\cX(o\ra o')} & & \rdTo^{\eta(o')} & \\
\cX(o) & & \dImplies_{\eta(o\ra o')} & & \cY(o') \\
& \rdTo_{\eta(o)} & & \ruTo_{\cY(o\ra o')} & \\
& & \cY(o). & &
\end{diagram}
\]
Note that in the 1-categorical setting, the commutativity of the above
diagram is precisely the condition that $\eta$ be a natural
transformation (the corresponding commutativity conditions on $\eta$
to be a natural transformation of 2-functors are too complicated to
write down here).

In particular, for every morphism $S:\emptyset \ra \emptyset$ from the
empty 1-manifold to itself (in $\cC$) we get a natural transformation
\[
\begin{diagram}[height=1.3em,width=3em,labelstyle=\scriptstyle]
& & \pt & & \\
& \ruTo^{\cX(S)} & & \rdTo^{\eta(\emptyset)} & \\
\pt & & \dImplies_{\eta(S)} & & \pt \\
& \rdTo_{\eta(\emptyset)} & & \ruTo_{\cY(S)} & \\
& & \pt. & &
\end{diagram}
\]
Since morphisms in $\cD$ between a point and itself are given by
$\D(\pt)$, and it is reasonable to expect that $\eta(\emptyset) =
\cO_\pt$, it follows that $\eta(S)$ should be thought of as a morphism
$\cX(S)\ra \cY(S)$ (in $\D(\pt)$, i.e., a morphism of graded vector
spaces).  In particular, taking $S=T^2$, it follows that any such
natural transformation will induce a map of graded vector spaces
\[ HH_*(X) = \cX(T^2) \ra \cY(T^2) = HH_*(Y). \]

\subsection{}
To summarize the above discussion, the main conjecture that we make is
that every map $\cE:\cX(S^1) \ra \cY(S^1)$ should lift to a natural
transformation of 2-functors $\eta_\cE:\cX \Rightarrow \cY$ (in an
appropriate sense), which in turn should induce a map of graded vector
spaces $\cX(S) \ra \cY(S)$ for every closed surface $S$, and in
particular the existence of the map
\[ (\FMXY^\cE)_* : HH_*(X) \ra HH_*(Y) \]
should be regarded as conjectural evidence for the existence of such a
lifting.

\subsection{}
As a quick check of this conjecture, note that it also predicts that
to an object $\cE\in \D(X\times Y)$ we should be able to associate in
a natural way a map $H^*(X, \cO_X) \ra H^*(Y, \cO_Y)$ (this
corresponds to taking $S=S^2$ instead of $T^2$ above): we believe the
map should be given by
\[
\begin{diagram}[height=2em,width=2em,labelstyle=\scriptstyle]
H^*(X, \cO_X) = \Hom_X(\cO_X, \cO_X)&\, & \rTo^{\FMXY^\cE} &
\Hom_Y(\FMXY^\cE(\cO_X), \FMXY^\cE(\cO_X)) & \rTo \\
&\, &  \rTo^{\Tr} & \Hom_Y(\cO_Y, \cO_Y) = H^*(Y, \cO_Y^*).&
\end{diagram}
\]
A particular consequence of this conjecture would thus be that the
numbers $h^{i,0}(X)$ should naturally be derived category invariants.

\subsection{}
The problem with making these conjectures precise (and attempting to
prove them) is the fact that it is obvious that one needs to change
the category $\cC$ and the construction of the functor $\cX$
associated to a space $X$ so as to incorporate the fact that $S_X$ is
nontrivial unless $X$ is a point.  For example, note that in our
association of functors to cobordisms, the left and right adjoints of
such functors only differ by a number of Serre functors (e.g.,
$\Delta^*$ is the left adjoint of $\Delta_*$, while the right adjoint
is isomorphic to $S_X^{-1} \circ \Delta^*$, for $\Delta:X\ra X\times
X$ the diagonal embedding).  The category $\cC$ needs to be changed in
such a way as to break the symmetry given by the fact that a
3-dimensional cobordism $S_1\Leftrightarrow S_2$ can be read as either
$S_1\Rightarrow S_2$ or $S_2\Rightarrow S_1$, for example by labeling
2-dimensional cobordisms by integers, and letting 3-cobordisms
``flow'' from the larger integer to the smaller.  Then to a labeled
2-cobordism we would still associate an object of $D(X\times\cdots
\times X)$ which would be supported on the small diagonal, but given
by $S_X^{\otimes n}$ for some appropriate $n$.

\subsection{}
\begin{figure}
\begin{center}
\includegraphics*{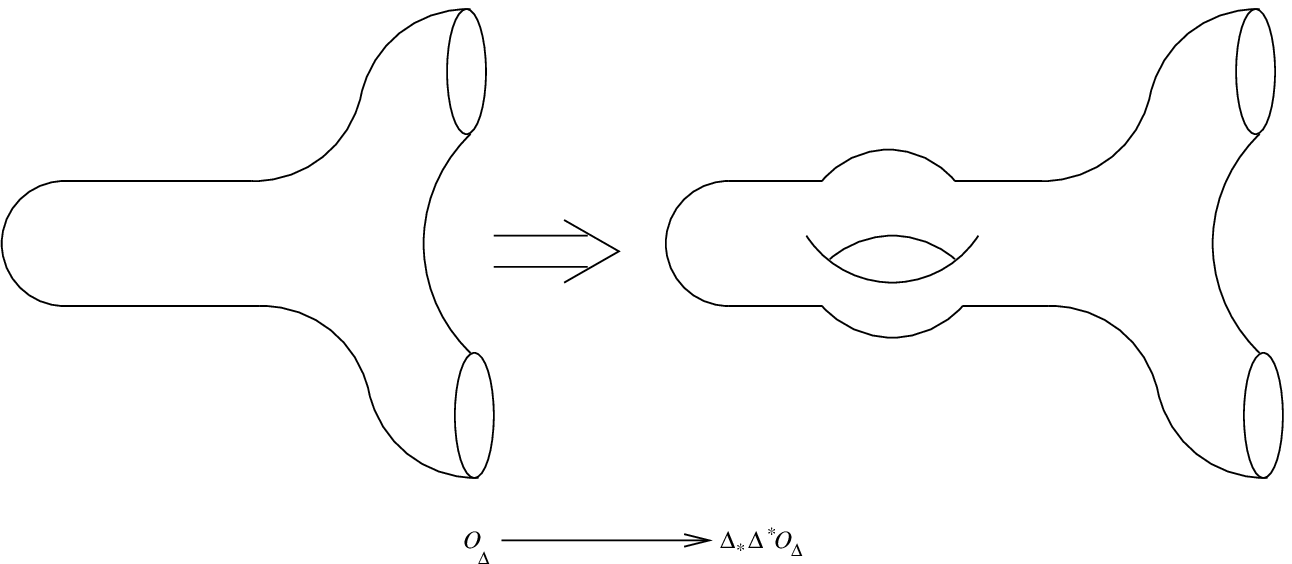}
\caption{The Atiyah class}
\label{fig:fig3}
\end{center}
\end{figure}

We conclude with a remark about the relevance of the Atiyah class in
the context of TQFT's.  This connects well with Kapranov's approach to
the Rozansky-Witten TQFT~\cite{Kap}.  The point is that we shall see
in~\cite{CalHH2} that the Atiyah class can be seen (via the HKR
isomorphism) to be nothing else but the unit
\[ \cO_\Delta \ra \Delta_*\Delta^* \cO_\Delta \]
of the adjunction $\Delta^*\adjoint\Delta_*$.  Pictorially this
corresponds to the cobordism depicted in Figure~\ref{fig:fig3},
induced by a standard surgery on a pair of pants, which is one of the
two fundamental building blocks of cobordisms of surfaces.

\section{DG-categories versus derived categories}
\label{app:categ2}

As mentioned in the introduction, working with derived categories as
triangulated categories raises several technical problems.  We discuss
these problems in this appendix, and we point out a possible way of
solving them.  The main difficulty arises from the fact that the
functor $\Phi$ that we introduced in Appendix~\ref{app:categ1} is not
fully faithful, even if restrict it to map into a smaller category
than $\gCat$, as the following example shows:

\begin{example}
\label{example:counter}
Let $X=Y=$ an elliptic curve over $\C$, let $\cF=\cO_\Delta$, the
structure sheaf of the diagonal in $X\times X$, and let
$\cG=\cO_\Delta[2]$.  Then $\Phi^\cF = \Id$, $\Phi^{\cG}$ is the
translation by 2 functor.  It is a straightforward calculation to see
that
\[ \Hom_{X\times Y}(\cF, \cG) = \Ext^2_{X\times Y}(\cO_\Delta,
\cO_\Delta) = \C, \] 
but the fact that $\gCoh(X)$ has cohomological dimension 1 implies that
any natural transformation between the identity functor and the
translation by 2 functor on $X$ must be zero.  (One has to use the
fact that every complex on $X$ is quasi-isomorphic to the direct sum
of its cohomologies.)  Therefore the functor $\Phi$ can not be
faithful in general.
\end{example}

\subsection{}
The trouble with the functor $\Phi$ as defined is that for two spaces
$X$ and $Y$, it is trying to be an equivalence between a triangulated
category ($\Hom(X, Y)$) and a category which has no obvious
triangulated structure (exact functors between $\D(X)$ and $\D(Y)$).
There is another category of functors which is naturally triangulated
and which we could use in our situation, namely
\[ \bExFun(X,Y) = H^0\Prex(\PreD(X), \PreD(Y)), \]
where $\Prex$ is the DG-category of preexact functors between the
DG-enhancements $\PreD(X)$, $\PreD(Y)$ of $\D(X)$, $\D(Y)$, in the
sense of Bondal and Kapranov~\cite{BonKap}.  The functor $\Phi$ can be
defined as before.  It seems reasonable to expect that $\bExFun(X, Y)$
is independent of the choice of enhancement, carries a Serre functor
which on objects is given by 
\[ F \mapsto S_Y\circ F\circ S_X, \]
and $\Phi$ is fully faithful.  As pointed out in
Appendix~\ref{app:categ1}, if this is correct we can think of objects
in $\D(X\times Y)$ as functors, and of morphisms between objects as
natural transformations, in the above sense.

\subsection{}
As an example of this kind of reasoning, the following is a more
intuitive ``definition'' of Hochschild homology and cohomology (see
also~(\ref{subsec:defhh})):
\begin{qdefinition}
\label{def:hhnat}
Define
\begin{align*}
HH^i(X) & = \Hom_{\bExFun(X,X)}(1_X, [i]), 
\intertext{where $1_X$ and $[1]$ are the identity and translation
functors on the pretriangulated category $\PreD(X)$~(\cite[Section
3]{BonKap}).  Similarly, define} 
HH_i(X) & = \Hom_{\bExFun(X,X)}(S_X^{-1}\circ [i], 1_X),
\end{align*}
where $S_X$ is the Serre functor of $\D(X)$ (which lifts to a functor
on $\PreD(X)$).
\end{qdefinition}

Again, cohomology is naturally a ring (with multiplication given by
composition of natural transformations) and homology is a module over
cohomology. 

\subsection{}
\label{subsec:mukproddual}
The isomorphism $\tau$ in the definition of the Mukai pairing has an
obvious categorical interpretation: it is just the isomorphism $\tau$
in~(\ref{subsec:dualiso}).  The Mukai pairing on $HH_*(X)$ can now be
reinterpreted by noting that the Serre dual of $HH_i(X)$ (with respect
to the Serre functor on $\bExFun(X,X)$) is the space
\[ \Hom_{\bExFun(X,X)}(\Id, S_X\circ [i]), \] 
which by~(\ref{subsec:dualiso}) is isomorphic to $HH_{-i}(X)$ via
$\tau$.  The same definition of the Mukai pairing can now be written
in the categorical context.

\bigskip \noindent
\small\textsc{Department of Mathematics, \\
University of Philadelphia, \\
Philadelphia, PA 19104-6395, USA} \\
{\em e-mail: }{\tt andreic@math.upenn.edu}

\end{document}